\documentclass[12pt,twoside]{article}

\usepackage{times}
\usepackage{amsmath,amssymb}
\usepackage{color}

\allowdisplaybreaks

\newcommand{\ef}{ \hfill $ \blacksquare $ \vskip 3mm}
\newcommand{\be}{\begin{equation}}
\newcommand{\ee}{\end{equation}}
\newcommand{\bea}{\begin{eqnarray}}
\newcommand{\eea}{\end{eqnarray}}

\newcommand{\bR}{{\mathbb R}}

\newcommand{\ov}{\overline{v}}
\def\nn{\nonumber}

\def\p{\partial}
\def\ve{\varepsilon}
\def\la{\lambda}

\def\al{\alpha}

\def\q{\quad}
\def\th{\theta}
\def\g{\gamma}
\def\G{\Gamma}
\def\dl{\delta}
\def\f{\frac}

\def\si{\sigma}
\def\b{\beta}

\def\q{\qquad}

\def\b{\beta}
\def\dl{\delta}
\def\Dl{\Delta}
\def\i{\infty}

\def\pa{\|}

\def\supp{\text{supp }}

\def\p{\partial}
\def\f{\frac}
\def\na{\nabla}
\def\al{\alpha}

\def\O{\Omega}

\def\q{\qquad}
\def\s{\sqrt}
 \allowdisplaybreaks

\begin{document}
 \footskip=0pt
 \footnotesep=2pt
\let\oldsection\section
\renewcommand\section{\setcounter{equation}{0}\oldsection}
\renewcommand\thesection{\arabic{section}}
\renewcommand\theequation{\thesection.\arabic{equation}}
\newtheorem{claim}{\noindent Claim}[section]
\newtheorem{theorem}{\noindent Theorem}[section]
\newtheorem{lemma}{\noindent Lemma}[section]
\newtheorem{proposition}{\noindent Proposition}[section]
\newtheorem{definition}{\noindent Definition }[section]
\newtheorem{remark}{\noindent Remark}[section]
\newtheorem{corollary}{\noindent Corollary}[section]
\newtheorem{example}{\noindent Example}[section]
\title{Regularity of Solution to Axis-symmetric Navier-Stokes Equations with a Slightly Supercritical Condition}

\author{Xinghong Pan
\footnote{E-mail:math.scrat@gmail.com.}\vspace{0.5cm}\\
  October 23, 2014\\
\small (Department of Mathematics and
IMS, Nanjing University, Nanjing 210093, China.)\\
\vspace{0.5cm}
}

\date{}
\maketitle
% \vskip 0.2in

\centerline {\bf Abstract} \vskip 0.3 true cm

Consider an axis-symmetric suitable weak solution of 3D incompressible Navier-Stokes equation with nontrivial swirl, $v=v_re_r+v_{\th}e_{\th}+v_ze_z$. Let $z$ denote the axis of symmetry and $r$ be the distance to the $z$-axis. If the solution satisfies a slightly supercritical assumption ( that is, $|v|\leq C\frac{(\ln|\ln r|)^\al}{r}$ for $\al$ $\in$ $[0,0.028]$ when $r$ is small ), then we prove that $v$ is  regular. This extends the results in \cite{Chen02},\cite{Koch01},\cite{Lei01} where regularities under critical assumptions, such as $|v|\leq \f{C}{r}$,  were proven.

\vskip 0.3 true cm

\vskip 0.3 true cm


\section{Introduction}
The incompressible Navier-Stokes equations in \textit{cartesian coordinates} are given by
\[
\p_t v+(v\cdot\nabla) v+\nabla p=\Delta v, \quad \nabla\cdot v=0,
\]
where $v$ is the velocity field and $p$ is the pressure. We consider the axis-symmetric solution of the equations. That means, in \textit{cylindrical coordinates} $r,\th,z$ with $x=(x_1,x_2,x_3)=(r\cos\th,r\sin\th,z)$, the solution is of this form
\[
v=v_re_r+v_{\th}e_{\th}+v_ze_z,
\]
where the basis vectors $e_r,e_\th,e_z$ are
\[
e_r=(\frac{x_1}{r},\frac{x_2}{r},0),\quad e_\th=(-\frac{x_2}{r},\frac{x_1}{r},0),\quad e_z=(0,0,1),
\]
and the components $v_r$,$v_\th$,$v_z$ do not depend on $\th$.

Recall $v_r,v_\th,v_z$ satisfy
\begin{equation}
\left\{
\begin{aligned}
&\p_t v_r+(b\cdot\nabla)v_r -\frac{(v_\th)^2}{r}+\p_r p=(\Delta-\frac{1}{r^2})v_r; \\
&\p_t v_\th+(b\cdot\nabla) v_\th+\frac{v_\th v_r}{r}=(\Delta-\frac{1}{r^2})v_\th ; \\
&\p_t v_z+(v\cdot\nabla)v_z+\p_z p=\Delta v_z ;                                    \\
&b=v_re_r+v_ze_z,\q \nabla\cdot b=\p_rv_r+\frac{v_r}{r}+\p_zv_z=0.
\end{aligned}
\right.\label{e1.1}
\end{equation}

In this paper we study the axis-symmetric Navier-Stokes equations under a slightly supercritical assumption on the drift term $b$. To be precise, we consider $b$ such that
\begin{equation}
|b|=\s{v^2_r+v^2_z}\leq\left\{
\begin{aligned}
&\frac{(\ln|\ln\frac{r}{3}|)^\al}{r}       & {\rm if}\ r \leq 1;\\
&    \frac{C}{r}                                & {\rm if}\ r>1.  \\
\end{aligned}
\right.\label{e1.2}
\end{equation}
Here $\al \in [0,0.028]$ is any fixed constant. Later we will see how $0.028$ is obtained.

Recall that the quantity $\G=rv_\th$ satisfies
\be
\p_t \G +(b\cdot\nabla)\G-\Delta\G +\frac{2}{r}\p_r \G=0.\label{e1.30}
\ee

Our assumption on $b$ is closely related to a counterexample in \cite{Seregin01}. In \cite{Seregin01}, the authors consider elliptic equation of this form
\be
-\Delta u+(b\cdot\nabla)u=0.\label{e1.40}
\ee
They construct a counterexample to state that \eqref{e1.40} does not have Liouville theorem when the divergence-free vector field $b$ satisfies $|b|\leq \f{\ln\ln|x|}{|x|}$ for large $|x|$. Morever, H\"{o}lder continuity, as well as Harnack inequality, to solutions of \eqref{e1.40} are also not to be expected. So under the assumption of \eqref{e1.2}, we do not expect a H\"{o}lder continuity to solutions of \eqref{e1.30} even if the exponent $\al$ is small.

Therefore, under the current techniques, the $''$lnln$''$ supercritical assumption on $b$ seems to be the best that one can expect for some continuity results(weaker than the h\"{o}lder continuity) to solutions of \eqref{e1.30} which can be used to prove the regularity of solutions to \eqref{e1.1}.

Our main result is the following.
\begin{theorem}\label{t1.1}
Let $(v,p)$ be a suitable weak solution of the axis-symmetric Navier-Stokes equation \eqref{e1.1} in $\bR^3\times [-1,0)$. Assume that $b$ satisfies \eqref{e1.2} and $\sup\limits_{x\in\bR^3}|\G(\cdot,-1)|< +\i$. Then we have
\[
\sup\limits_{(x,t)\in\bR^3\times[-1,0)}|v|< +\i.
\]
\end{theorem}
Here and throughout the paper, we will use $c$ and $C$ to denote a generic constant. It may be different from line to line. Also we use $A\lesssim B$ to denote $A\leq CB$.
\begin{remark}\label{r1.1}
 We note that $\G$ satisfies the equation \eqref{e1.30} which enjoys the maximal principle. So the assumption $\sup\limits_{x\in\bR^3}|\G(\cdot,-1)|< +\i$ can assure that $|\G|\leq C$ for all $t\in [-1,0)$ for some positive constant $C$.

 Readers can refer to \cite{Chen01} for the definition of suitable weak solutions.
\end{remark}

Recall the natural scaling of Navier-Stokes equations: If $(v,p)$ is a solution of equations \eqref{e1.1}, then for any $\la>0$, the following rescaled pair is also a solution:
\be
v^\la(x,t)=\la v(\la x,\la^2t),\q p^\la(x,t)=\la^2p(\la x,\la^2t).  \nn
\ee

Denote $b^\la=v^\la_re_r+v^\la_ze_z$, we say our assumption \eqref{e1.2} is supercritical which means, for a fixed point $x_0=(r_0\cos\th_0,r_0\sin\th_0,z_0)$, $b^\la(x_0)$ satisfies a bound
\be
|b^\la(x_0)|\leq \frac{(\ln|\ln\la\frac{r_0}{3}|)^\al}{r_0}. \nn
\ee
When $\la\rightarrow 0$, the bound goes to infinity. That is, when one zooms in at a point , the bound on the drift term becomes worse, so the regularity of our solution must be handled more carefully.

Global in-time regularity of the solution to the axis-symmetric Navier-Stokes equations is still open. Under the no swirl assumption, $v^\th=0$ , Ladyzhenskaya\cite{Ladyzenskaja01} and Ukhovskii-Yudovich \cite{Ukhovskii01} independently proved that weak solutions are regular for all time. When the swirl $v^\th$ is non-trivial, recently, some efforts and progress have been made on the regularity of the axis-symmetric solutions. In \cite{Chen02}, Chen-Strain-Yau-Tsai proved that the suitable weak solutions are regular if the solution satisfies $|v|\leq C/r<\i$. Their method is based on Nash\cite{Nash01},Moser\cite{Moser01} and De Giorgi\cite{De Giorgi01}. Also, Koch-Nadirashvili-Seregin-Sverak in \cite{Koch01} proved the same result using a Liouville theorem and scaling-invariant property. Lei and Zhang in \cite{Lei01} proved regularity of the solution under a more general assumption on the drift term $b$ where $b\in L^\i\left([-1,0),BMO^{-1}\right)$.

It seems that their assumptions on $b$ are critical( for a fixed point, after scaling, the bound on $b$ is invariant ). So using a standard linear estimate, they can prove the H\"{o}lder continuity of $\G$ from equation \eqref{e1.30} which breaks the scaling-invariant bound of the angular component $v_\th$, making the bound on $b$ to a subcritical one. This is very important in proving regularity of the solution $v$. But under our supercritical assumption \eqref{e1.2}, only a logarithmic modulus of continuity, rather than  the H\"{o}lder continuity, can be obtained which indicates, near $r=0$, the $L^\i$ norm of $\G$ has a logarithmic decay with respect to $r$. Note that this also breaks the scaling-invariant bound of $v_\th$ and is enough to prove the regularity of $v$, but requiring more efforts and more complicated computation.

Our proof of Theorem1.1 is initially inspired by \cite{Chen02} and \cite{Lei01}. In the appendix of \cite{Chen02}, the authors give a time-independent bound to the axis-symmetric weak solution under the assumption $|v|\leq C/r$. Stimulated by their idea, we will give a similar proof to get the regularity of $v$ under the assumption \eqref{e1.2}. In the process, more detailed computation and careful handling will be needed , especially when we deal with the estimate to the fundamental solution of  \eqref{e1.4} due to the critical term $\f{2}{r}\p_r$. Our procedures of proof are as follows.

First, we will follow \cite{Lei01}, using Nash-Moser type method to prove continuity of $\G$ at $r=0$. It satisfies a $\log$ decay near $r=0$, that is
\be
|\G|\leq C|\ln {\frac{r}{3}}|^{-c_0}  \q {\rm when}\  r\leq1, \label{e1.50}
\ee
for some small positive $c_0$. See Theorem 1.2. This estimate breaks the scaling-invariant bound of $v_\th$.

Next we explore the relationship between $v_\th$ and $w_\th$, the angular component of the vorticity $w=\nabla \times v$.\\ Here
\[
w(x,t)=w_re_r+w_\th e_\th+w_ze_z
\]
and
\[
w_r=-\p_z v_\th, \ w_\th=\p_z v_r-\p_r v_z,\  w_z=(\p_r+\frac{1}{r})v_\th,
\]
where $w_\th$ satisfies
\begin{equation}
\left[\p_t+b\cdot\nabla-\Delta-\frac{v_r}{r}\right]w_\th-\p_z\frac{(v_\th)^2}{r}+\frac{w_\th}{r^2}=0. \nn
\end{equation}
Let $\O=\f{w_\th}{r}$, then $\O$ satisfies
\be
(\p_t-L)\O=r^{-2}\p_z(v_\th)^2, \q L=\Dl+\f{2}{r}\p_r-b\cdot\na.    \label{e1.4}
\ee

Combining \eqref{e1.50} and \eqref{e1.4}, we can get an estimate of $w_\th$. Handling of \eqref{e1.4} involves an estimate to the fundamental solution of $(\p_t-L)u=0$ which will be described in Theorem 1.3.

Recall that $b$ satisfies the vector identity
\be
-\Delta b={\rm curlcurl}b-\na{\rm div}b,  \nn
\ee
and note that
\be
{\rm div}b=0,\q {\rm curl}b=w_\th e_\th,  \nn
\ee
then we have
\begin{equation}
-\Delta b=\na\times(w_\th e_\th).   \nn
\end{equation}

At last, using the regularity theory of elliptic equations, we can get the boundedness of $b$. This will prove the regularity of our solution $v$.
\begin{theorem}
For the divergence-free drift term $b$, define a zero-dimensional integral norm
\begin{equation}
E_R(b)\triangleq \sup\limits_{-R^2\leq t\leq 0} \left\{\frac{1}{R^{3-p}}\int_{B_R}|b|^p dx\right\}^{\frac{1}{p}},\label{e1.5}
\end{equation}
where $\frac{5}{3}< p \leq 2$ and $R\leq 1$, $B_R$ is the ball of radius $R$ centered at $x=0$.If
\begin{equation}
E_R(b)\leq C\left(\ln|\ln\frac{R}{3}|\right)^{\frac{3p-5}{77p-120}(1-\b)}\q \forall R\in(0,1],\label{e1.6}
\end{equation}
for any $\b>0$, then the weak solution of  \eqref{e1.30} is  continuous at $(0,0)$ and it has a log decay near $r=0$. That is, there exists some positive $c_0$, such that
\[
|\G|\leq  C|\ln{\frac{r}{3}}|^{-c_0}   \q{\rm for}\ r\leq 1.
\]
\end{theorem}
\begin{remark}\label{r1.2}
We note that when $|b|\leq \frac{(\ln|\ln\frac{r}{3}|)^\al}{r}$,$\al\in[0,0.028]$, there exists a $p_0 \in (\frac{5}{3},2]$ and $\b_0$ small such that the assumption \eqref{e1.6} is satisfied. So Theorem1.2 infers the log continuity of $\G$. We compute it as follows,
\bea
E^p_R(b)&\leq& \frac{1}{R^{3-p}}\int_{B_R}\left[\frac{(\ln|\ln\frac{r}{3}|)^\al}{r}\right]^p dx \nn    \\
      &\leq& \frac{2\pi}{R^{3-p}}\int^R_0\int^R_0 \left[\frac{(\ln|\ln\frac{r}{3}|)^\al}{r}\right]^pr dr dz \nn \\
      &\leq& \frac{C}{R^{2-p}}\int^R_0 (\ln|\ln\frac{r}{3}|)^{\al p} r^{1-p} dr \nn   \\
      &\leq& \frac{C}{R^{2-p}}\int^\infty_{\frac{1}{R}} (\ln\ln3s)^{\al p} s^{p-3} ds \q\q (\rm{let}\  s=1/r)   \nn \\
      &\leq&\frac{C}{R^{2-p}}\left(\int^{\frac{1}{3}e^{e^{\left(\ln\ln\frac{3}{R}\right)^{1+\frac{\ve}{\al}}}}}_{\frac{1}{R}}+\int^\infty_
      {\frac{1}{3}e^{e^{\left(\ln\ln\frac{3}{R}\right)^{1+\frac{\ve}{\al}}}}}\right) (\ln\ln3s)^{\al p} s^{p-3} ds\q \rm{for\  any}\ \ve>0    \nn   \\                                                     &\leq& \frac{C}{R^{2-p}}\left\{\int^{\frac{1}{3}e^{e^{\left(\ln\ln\frac{3}{R}\right)^{1+\frac{\ve}{\al}}}}}_{\frac{1}{R}} (\ln\ln3s)^{\al p}s^{p-3} ds\right.\nn\\ &&\left.+\int^\infty_{\frac{1}{3}e^{e^{\left(\ln\ln\frac{3}{R}\right)^{1+\frac{\ve}{\al}}}}} \left(\frac{\ln\ln3s}{3s}\right)^{\al p}\left(3s\right)^{p-3+\al p} ds \right\}. \nn                                                                        \eea
Here $\ln\ln3s$ is a monotone-increasing function while $\frac{\ln\ln3s}{3s}$ a monotone-decreasing function , so
\begin{flalign}
E^p_R(b)\leq&\frac{C}{R^{2-p}}\left\{ \left(\ln\ln\frac{3}{R}\right)^{\al p+\ve p} \int^{\frac{1}{3}e^{e^{\left(\ln\ln\frac{3}{R}\right)^{1+\frac{\ve}{\al}}}}}_{\frac{1}{R}} s^{p-3} ds\right.&\nn\\
 &\left.+\frac{\left(\ln\ln\frac{3}{R}\right)^{\al p+\ve p}}{\left(e^{e^{\left(\ln\ln\frac{3}{R}\right)^{1+\frac{\ve}{\al}}}}\right)^{\al p}}\int^\infty_{\frac{1}{3}e^{e^{\left(\ln\ln\frac{3}{R}\right)^{1+\frac{\ve}{\al}}}}} \left(3s\right)^{p-3+\al p} ds \right\},& \nn
\end{flalign}
we need choose a $p \in (\frac{5}{3},2] $ such that $p-3+\al p <-1$,that is $\al<\f{2}{p}-1$.

 For such a $p$
\bea
E^p_R(b)&\leq& \frac{C}{R^{2-p}}\left\{ \left(\ln\ln\frac{3}{R}\right)^{\al p+\ve p} \frac{1}{2-p}\left(\frac{1}{R}\right)^{p-2}\right.      \nn  \\
      &&\left.+\frac{1}{2-\al p-p}\frac{\left(\ln\ln\frac{3}{R}\right)^{\al p+\ve p}}{\left(e^{e^{\left(\ln\ln\frac{3}{R}\right)^{1+\frac{\ve}{\al}}}}\right)^{\al p}}\left(e^{e^{\left(\ln\ln\frac{3}{R}\right)^{1+\frac{\ve}{\al}}}}\right)^{p-2+\al p}\right\},       \nn
\eea
here $e^{e^{\left(\ln\ln\frac{3}{R}\right)^{1+\frac{\ve}{\al}}}}\geq\frac{3}{R}$. So
\bea
E^p_R(b)&\leq& \frac{C}{R^{2-p}} \left(\ln\ln\frac{3}{R}\right)^{\al p+\ve p} \left(\frac{1}{R}\right)^{p-2}               \nn  \\
      &\leq& C\left(\ln\ln\frac{3}{R}\right)^{\al p+\ve p}.         \nn
\eea
Let $\ve$ and $\b_0$ are sufficiently small to make $\al+\ve\leq \frac{3p-5}{77p-120}(1-\b_0)$.
This means, $\al<\frac{3p-5}{77p-120}$. So
\be
\al<\min\left\{\frac{3p-5}{77p-120},\f{2}{p}-1\right\}.  \nn
\ee
Let $f(p)=\min\left\{\frac{3p-5}{77p-120},\f{2}{p}-1\right\} $, $f(p_0)=\max\limits_{\f{5}{3}<p\leq2}f(p)$. We compute that
\be
p_0=\f{279+\s{1041}}{160}\approx 1.945.     \nn
\ee
We choose such $p_0$ to ensure $\al<\min\left\{\frac{3p_0-5}{77p_0-120},\f{2}{p_0}-1\right\}\approx0.028$. This is nearly the maximum value we can choose for $\al$.
\end{remark}

The next theorem gives an upper bound estimate to the fundamental solution of equation
\be
\p_t u=(\Dl+\f{2}{r}\p_r-b\cdot\na)u, \q \na\cdot b=0\label{1.9}
\ee
under certain bound for $b$, which will be useful in the proof of Theorem 1.1. Due to the term $\f{2}{r}\p_r$, the result is not covered by the standard theory.

Before stating the theorem, we give the definition of fundamental solutions to \eqref{1.9}.
\begin{definition}
Let $Q=\{(x,t)|x\in \bR^3, t>s\}$, we say $0\leq p(x,t;y,s)\in L_{\rm loc}(Q)\bigcap C^2(\overline{Q}\backslash (y,s))$ is a fundamental solution of \eqref{1.9} in $Q$ if it satisfies\\
1. for any $\psi(x,t)\in C^{\i}_0(Q)$ and $\psi|_{r=0}=0$,
\be
\int^{+\i}_s\int_{\bR^3}p(x,t;y,s)(\psi_t+\Delta\psi+b\cdot\na\psi-\f{2}{r}\p_r\psi)dxdt=0.  \nn
\ee
2. for any $\phi(x)\in C^{\i}_0(\bR^3)$,
\be
\lim\limits_{t\rightarrow s}\int_{\bR^3}p(x,t;y,s)\phi(y)dy=\phi(x). \nn
\ee
3. Let $y=(y_1,y_2,y_3)$ and denote $y'=(y_1,y_2,0)$, we require
\be
p(x,t;y,s)|_{|y'|=0}=0. \nn
\ee

 This third condition marks an important difference with the standard theory where fundamental solutions are positive everywhere. Our choice of this fundamental solution coincides with some quantities in the axis-symmetric Navier-Stokes equations, such as $\G$, $w_\th$.
\end{definition}
\begin{remark}
Due to our assumption $p(x,t;y,s)\in C^2(\overline{Q}\backslash (y,s))$, $p(x,t;y,s)$ satisfies \eqref{1.9} in classical sense except at point $(y,s)$.
\end{remark}
\begin{theorem}
Let $p(x,t;y,s)$ be a fundamental solution of \eqref{1.9} and the divergence-free smooth vector function $b(x,t)$ satisfies $|b|\leq C_0+\frac{1}{r}$. Then we have
\be
p(x,t;y,s)\leq C(t-s)^{-3/2}\exp\left\{-C_1\frac{|x-y|^2}{t-s}\left(1-C_0\f{t-s}{|x-y|}\right)^2_{+}\right\} \label{1.10}
\ee
for some positive constants $C,C_1$. Moreover,
\be
\int_{\bR^3}p(x,t;y,s)dx\leq 1,\q \int_{\bR^3}p(x,t;y,s)dy= 1.  \nn
\ee
\end{theorem}
\begin{remark}
The idea of proving Theorem 1.3 is based on Theorem 5 of \cite{Carlen01}, but due to the term $\f{2}{r}\p_r$, the proof will be more complicated. In \cite{Carlen01}, the authors consider the equation
\be
\p_t u=\Dl u-b\cdot\na u.  \label{e1.100}
\ee
In their proof, the Davies-type exponent $r(t)$ can map from $[0,T]$ to $[1,\i)$ and with the help of the logarithmic Sobolev inequality, a $L^1\rightarrow L^\i$ estimate to the solution of \eqref{e1.100} can be obtained. But for \eqref{1.9}, we must deal with a singular term $\f{2}{r}\p_r$ which will create some difficulties when using their method to estimate the solution. However, stimulated by Fabes-Stroock\cite{Fabes01} and Davies \cite{Davies01}, we use a dual technique from harmonic analysis to overcome this difficulty. We proceed as follows.

First, we choose $r(t): [0,T]\rightarrow[2,\i)$ to get a $L^2\rightarrow L^\i$ estimate of the fundamental solution $p(x,t;y,s)$, then the same estimate can be applied to the adjoint $p^*(x,t;y,s)$ of $p(x,t;y,s)$. By duality, we get $L^1\rightarrow L^\i$ estimate of $p(x,t;y,s)$. This will prove our Theorem 1.3.

Estimates to the kernels of parabolic equations have had a long history especially when the drift $b$ is a divergence-free singular term.
Under different assumptions on $b$ , Osada H.\cite{Osada01}, Liskevich-Zhang\cite{Lis01},Zhang Qi S.\cite{Zhang01} give bounds for the fundamental solution of \eqref{e1.100}. Readers can refer to their papers and their References for more information. Here we add a singular term $\f{2}{r}\p_r$ in the equation and give an upper bound to the fundamental solution. We hope our estimate can not only be applied to the axis-symmtric Navier-Stokes equations, but also to other related incompressible fluid fields.
\end{remark}

We now recall some regularity results on the axis-symmetric Navier-Stokes equations. In the presence of swirl, from the partial regularity theory of \cite{Caffarelli01}, any singular points of the axis-symmetric suitable weak solution of \eqref{e1.1} can only lie on the $z$ axis.In \cite{Buiker01}, Burke-Zhang give a priori bounds for the vorticity of axially symmetric solutions which indicates that the result of \cite{Caffarelli01} can be applied to a large class of weak solutions. Chan-Vasseur in \cite{Vasseur01} give a logarithmically improved Serrin criterion for global regularity to solutions of Navier-Stokes equations. See also an extension in Zhou-Lei \cite{Zhou01}. Neustupa and Pokorny[6] proved certain regularity of one component(either $v^\th$ or $v^r$) imply regularity of the other components of the solutions. Chae-Lee[11] proved regularity assuming a zero-dimensional integral norm on $w^\th$: $w^\th\in L^s_tL^q_x$ with $3/q +2/s=2$. Also regularity results come from the work of Jiu-Xin [12] under the assumption that another zero-dimensional scaled norms $\int_{Q_R}(R^{-1}|w^\th|^2+R^{-3}|v^\th|^2) dz$ is sufficiently small for $R>0$ is small enough. On the other hand, Lei-Zhang\cite{Lei02} give a structure of singularity of 3D axis-symmetric equation near maximum point. Tian-Xin\cite{Tian01} constructed a family of singular axi-symmetric solutions with singular initial datas. Recently, Hou-Li\cite{Hou02} construct a special class of global smooth solutions. See also a recent extension: Hou-Lei-Li\cite{Hou01}.

The paper is organized as follows:In section 2, we establish a local maximum estimate using a Moser's iteration. Based on the local maximal estimate, In section 3 , we obtain the continuity of $\G$ and prove Theorem1.2 by Nash's method. In section 4,we prove Theorem1.1. The argument is based on \cite{Chen02}. In section 5,we give the proof of Theorem 1.3.
\section{Local Maximum Estimate}

In this section, Using Moser's iteration, we prove a local maximum estimate of $\G$ which will be used to obtain continuity of $\G$ in the next section. the main idea is to exploit the divergence-free property of $b(x,t)$ and a special cut-off function. We learn from Lei-Zhang\cite{Lei01} and \cite{Chen01} where the authors treated the term $\frac{2}{r}\p_r \G$ and $b\cdot\nabla \G$.

We first derive a parabolic De Giorgi type energy estimates of \eqref{e1.30}. Set $\frac{1}{2}\leq \si_2<\si_1\leq 1$ and consider a test function $\psi(y,s)=\phi(|y|)\eta(s)$ satisfying
\begin{equation}
\left\{
\begin{aligned}
&\supp\phi \subset B(\si_1), \phi=1 \ \textrm{in}\  B(\si_2), 0\leq\phi\leq 1;\\
&\supp\eta\subset(-(\si_1)^2,0], \eta=1 \ \textrm{in}\  (-(\si_2)^2,0], 0\leq\eta\leq1;  \\
& |\eta'|\lesssim \frac{1}{(\si_1-\si_2)^2}, \mid \frac{\nabla\phi}{\sqrt{\phi}}\mid \leq  \frac{1}{\si_1-\si_2}.              \\
\end{aligned}
\right.\label{e2.1}
\end{equation}

We will also use the following notations to denote our domains. Let $R>0$ and $R\in (0,1)$. We write $B_R=B(0,R)$ and
\[
P(R)=B_R\times(-R^2,0],\q P(R_1,R_2)=B_{R_1}/B_{R_2}\times(-R^2_1,0]\  \textrm{for} \ R_1>R_2.
\]

Consider the function $u=|\G|^q,\ q>1$ and the test function $\psi_R(y,s)=\phi_R(|y|)\eta_R(s)=\phi(\frac{y}{R})\eta(\frac{s}{R^2})$. Testing \eqref{e1.30} by $q|\G|^{2q-2}\G\psi^2_R$ gives
\begin{equation}
\frac{1}{2}\int\int (\p_s u^2+b\cdot\nabla u^2+\frac{2}{r}\p_r u^2)\psi^2_R=q\int\int\Delta \G |\G|^{2q-2}\G\psi^2_R. \label{e2.2}
\end{equation}
Using Cauchy-Schwartz's inequality and integration by parts, we compute
\bea
&&q\int\int\Delta \G |\G|^{2q-2}\G\psi^2_R     \nn  \\
&=&q\int\int\Delta |\G| |\G|^{2q-1}\psi^2_R    \nn  \\
&=&-q\int\int (2q-1)|\nabla\G|^2\G^{2q-2}\psi^2_R+|\nabla\G||\G|^{2q-1}\nabla\psi^2_R    \nn    \\
&=&-\int\int\frac{2q-1}{q}|\nabla\G^q|^2\psi^2_R+2\psi_R|\G^q|\nabla\psi_R\cdot|\nabla\G^q|  \nn \\
&=&-\int\int\frac{2q-1}{q}|\nabla u|^2\psi^2_R+2\psi_Ru\nabla\psi_R\cdot\nabla u  \nn     \\
&=&-\int\int(2-\frac{1}{q})|\nabla(u\psi_R)|^2-(2-\frac{2}{q})u\nabla \psi_R\cdot\nabla(u\psi_R)-\frac{1}{q}u^2|\nabla\psi_R|^2   \nn    \\
&\lesssim& -\int\int|\nabla(u\psi_R)|^2+\int\int u^2|\nabla \psi_R|^2 \label{e2.3}
\eea
and
\bea
\frac{1}{2} \int\int\p_s u^2 \psi^2_R=\frac{1}{2}\int_{B(\si_1R)}u^2(\cdot,t)\psi^2_R-\frac{1}{2}\int\int u^2\p_s\psi^2_R. \label{e2.4}
\eea
Moreover, by the fact that $\G=0$ on the axis $r=0$, we have
\bea
&&-\frac{1}{r}\int\int\p_r u^2\psi^2_R     \nn        \\
&=&-2\pi\int^{+\infty}_{-\infty}\int^{+\infty}_0\p_r u^2\psi^2_R dr dz  \nn    \\
&=& 2\pi\int^{+\infty}_{-\infty}\int^{+\infty}_0 u^2\p_r\psi^2_R dr dz   \nn   \\
&=& -2\pi\int^{+\infty}_{-\infty}\int^{+\infty}_0 \p_r(u^2\p_r\psi^2_R)r dr dz  \nn \\
&\lesssim& \int\int u^2|\p_r(\p_r \psi^2_R)|+u\p_r u\psi\p_r \psi   \nn   \\
&\lesssim& \int\int u^2(|\nabla\psi_R|^2+|\nabla^2\psi_R|)+u\p_r\psi_R(\p_r(u\psi_R)-u\p_r\psi_R) \nn   \\
&\lesssim& \int\int u^2(|\nabla\psi_R|^2+|\nabla^2\psi_R|)+\frac{1}{4}\int\int|\nabla(u\psi_R)|^2.   \label{e2.5}
\eea
Consequently, using \eqref{e2.2} and combining \eqref{e2.3},\eqref{e2.4} and \eqref{e2.5},we get
\bea
&&\int_{B(\si_1R)}u^2(\cdot,t)\psi^2_R+\int\int|\nabla(u\psi_R)|^2   \nn  \\
&\lesssim&\int\int u^2(|\nabla\psi_R|^2+|\nabla^2\psi_R|+|\p_s\psi^2_R|)-\frac{1}{2}\int\int b\cdot\nabla u^2\psi^2_R  \nn \\
&\lesssim&\frac{1}{(\si_1-\si_2)^2R^2}\int\int_{P(\si_1R)} u^2-\frac{1}{2}\int\int b\cdot\nabla u^2\psi^2_R.  \label{e2.6}
\eea
By the divergence-free property of the drift term $b$ and using \eqref{e2.1}, we have
\bea
&&-\frac{1}{2}\int\int b\cdot\nabla u^2\psi^2_R  \nn   \\
&=&\int\int u^2\psi_Rb\cdot\nabla\psi_R   \nn   \\
&=&\int\int(\psi_R u)^{3/2}u^{\frac{1}{2}}b\cdot\frac{\nabla \psi_R}{\sqrt{\psi_R}}   \nn    \\
&\lesssim&\frac{1}{(\si_1-\si_2)R}\int\int_{P(\si_1R,\si_2R)}|b|(\psi_R u)^{3/2}u^{\frac{1}{2}}   \nn     \\
&\lesssim&\frac{1}{(\si_1-\si_2)R}\int\left[\left(\int|b|^p dy\right)^{\frac{1}{p}}\left(\int(\psi_R u)^6dy\right)^{\frac{1}{4}}\left(\int u^{\frac{2p}{3p-4}}dy\right)^{\frac{3p-4}{4p}}\right] ds  \nn \\
&\lesssim&\frac{\parallel b\parallel_{L^\infty_tL^p_x(P(\si_1R))}}{(\si_1-\si_2)R}\int \left[\left(\int|\nabla(\psi_R u)|^2dy\right)^\frac{3}{4}\left(\int u^{\frac{2p}{3p-4}}dy\right)^{\frac{3p-4}{4p}}\right]ds                          \nn                     \\
&\lesssim&\frac{\parallel b\parallel_{L^\infty_tL^p_x(P(\si_1R))}}{(\si_1-\si_2)R}\left(\int\int|\nabla(\psi_R u)|^2dyds\right)^\frac{3}{4}\left(\int\int u^{\frac{2p}{3p-4}}dyds\right)^{\frac{3p-4}{4p}}\left(\int^{-(\si_1 R)^2}_{-(\si_2 R)^2}ds\right)^{\frac{2-p}{2p}}   \nn    \\
&\lesssim&\left[\frac{((\si_1-\si_2)^2R^2)^{\frac{2-p}{2p}}}{(\si_1-\si_2)R}\parallel b\parallel_{L^\infty_tL^p_x(P(\si_1R))}\right]^4\left(\int\int u^{\frac{2p}{3p-4}}\right)^{\frac{3p-4}{p}}+\frac{1}{2}\int\int|\nabla(\psi_R u)|^2             \nn           \\
&\lesssim&\left[\frac{1}{[(\si_1-\si_2)R]^{3-\frac{5}{p}}}E_R(b)\right]^4\left(\int\int u^{\frac{2p}{3p-4}}dyds\right)^{\frac{3p-4}{p}}+\frac{1}{2}\int\int|\nabla(\psi_R u)|^2dyds. \label{e2.7}
\eea
Combining \eqref{e2.6} and \eqref{e2.7},using the Cauchy-Schwartz inequality, we get
\bea
&&\int_{B(\si_1R)}u^2(\cdot,t)\psi^2_R+\int\int_{P(\si_1R)}|\nabla(u\psi_R)|^2      \nn    \\
&\lesssim&\frac{1}{(\si_1-\si_2)^2R^2}\int\int_{P(\si_1R)} u^2dyds+\frac{E^4_R(b)}{[(\si_1-\si_2)R]^{12-\frac{20}{p}}}\left(\int\int_{P(\si_1 R)} u^{\frac{2p}{3p-4}}dyds\right)^{\frac{3p-4}{p}}                                   \nn              \\
&\lesssim&\frac{1}{(\si_1-\si_2)^2R^2}\left(\int\int_{P(\si_1R)}u^{2\frac{p}{3p-4}}dyds\right)^{\frac{3p-4}{p}}\left(\int\int_{P(\si_1R,\si_2R)}dyds\right)^{\frac{4-2p}{p}}  \nn \\
&&+\frac{E^4_R(b)}{[(\si_1-\si_2)R]^{12-\frac{20}{p}}}\left(\int\int_{P(\si_1 R)} u^{\frac{2p}{3p-4}}dyds\right)^{\frac{3p-4}{p}}                  \nn \\
&\lesssim&\frac{1+E^4_R(b)}{[(\si_1-\si_2)R]^{12-\frac{20}{p}}}\left(\int\int_{P(\si_1 R)} u^{\frac{2p}{3p-4}}dyds\right)^{\frac{3p-4}{p}}.             \nn
\eea
At last,we get the estimate
\begin{flalign}
\hspace{10mm}
&\sup\limits_{-\si^2_1R^2\leq t\leq 0}\int_{B(\si_1R)}u^2(\cdot,t)\psi^2_R+\int\int_{P(\si_1R)}|\nabla(u\psi_R)|^2 &                      \nn     \\
\lesssim&\frac{(1+E_R(b))^4}{[(\si_1-\si_2)R]^{12-\frac{20}{p}}}\left(\int\int_{P(\si_1 R)} u^{\frac{2p}{3p-4}}dyds\right)^{\frac{3p-4}{p}}.&  \label{e2.8}
\end{flalign}

Our next step is to derive a mean value inequality based on \eqref{e2.8} using Moser's iteration.
\begin{lemma}
Suppose $u$ satisfies \eqref{e2.8} for $p\in (\frac{5}{3},2]$, then for $0<R\leq 1$, there is the estimate
\be
\sup\limits_{P(\frac{1}{2}R)} \G\lesssim (1+E_R(b))^{\frac{5p}{2(3p-5)}}\left(\int\int_{P(R)}\frac{1}{R^5}\G^2\right)^{\frac{1}{2}}. \label{e2.9}
\ee
\end{lemma}
\noindent{\bf Proof.}By H\"{o}lder inequality and Sobolev imbedding theorem, we have
\begin{flalign}
\hspace{25mm}
&\int\int_{P(\si_1R)}(u\psi_R)^{\frac{10}{3}} &                              \nn              \\
\lesssim& \int\left(\|u\psi_R(\cdot,s)\|^{\frac{4}{3}}_{L^2_{(B(\si_1R))}}\|\nabla(u\psi_R)\|^2_{L^2_{B(\si_1R)}}\right) ds & \nn   \\
\lesssim& \sup\limits_{-(\si_1R)^2\leq t\leq 0} \|u\psi_R(\cdot,t)\|^{\frac{4}{3}}_{L^2_{B(\si_1R)}}\|\nabla(u\psi_R)\|^2_{L^2_{P(\si_1R)}}.& \nn
\end{flalign}
Using \eqref{e2.1} and \eqref{e2.8},we get
\be
\int\int_{P(\si_2R)}u^{\frac{10}{3}} \lesssim \left\{\frac{(1+E_R(b))^4}{[(\si_1-\si_2)R]^{12-\frac{20}{p}}}\right\}^{\frac{5}{3}} \left(\int\int_{P(\si_1R)}u^{\frac{2p}{3p-4}}\right)^{\frac{5(3p-4)}{3p}}.  \nn
\ee
Remember $u=\G^q$, then we obtain
\bea
&&\int\int_{P(\si_2R)}\left(\G^{\frac{2p}{3p-4}q}\right)^{\frac{3p-4}{2p}\times\frac{10}{3}}    \nn      \\
&\lesssim& \left\{\frac{(1+E_R(b))^4}{[(\si_1-\si_2)R]^{12-\frac{20}{p}}}\right\}^{\frac{5}{3}}\left(\int\int_{P(\si_1R)}\G^{\frac{2p}{3p-4}q}\right)^{\frac{5(3p-4)}{3p}}. \label{e2.10}
\eea
For convenience of computation, we denote $\kappa=\f{p}{3p-4}$, then \eqref{e2.10} is
\bea
&&\int\int_{P(\si_2R)}\left(\G^{2\kappa q}\right)^{\frac{5}{3}\kappa^{-1}}dyds    \nn      \\
&\lesssim& \left\{\frac{(1+E_R(b))^4}{[(\si_1-\si_2)R]^{12-\frac{20}{p}}}\right\}^{\frac{5}{3}}\left(\int\int_{P(\si_1R)}\G^{2\kappa q}dyds\right)^{\frac{5}{3}\kappa^{-1}}. \label{e2.101}
\eea
For integer $j\geq 0$ and a constant $\si=\frac{1}{2}$, set $\si_2=\frac{1}{2}(1+\si^{j+1})$ and $\si_1=\frac{1}{2}(1+\si^j)$. Let $q=(\frac{5}{3}\kappa^{-1})^j$,then we get
\bea
&&\left(\int\int_{P(\frac{R}{2}(1+\si^{j+1}))} \G^{{2\kappa}(\frac{5}{3}\kappa^{-1})^{j+1}}dyds \right)^{\frac{1}{2\kappa}(\frac{3}{5}\kappa)^{j+1}}   \nn    \\
&\lesssim&\frac{(1+E_R(b))^{\frac{20}{3}\frac{1}{2\kappa}(\frac{3}{5}\kappa)^{j+1}}}{[\si^{(j+1)}R]^{\frac{5}{3}(12-\frac{20}{p})\frac{1}{2\kappa}(\frac{3}{5}\kappa)^{j+1}
          }}\times                    \nn                  \\
&&\left(\int\int_{P(\frac{R}{2}(1+\si^j))} \G^{{2\kappa}(\frac{5}{3}\kappa^{-1})^j}dyds \right)^{\frac{1}{2\kappa}(\frac{3}{5}\kappa)^j}. \nn
\eea
By iterating $j$, the above inequality gives
\bea
&&\left(\int\int_{P(\frac{R}{2}(1+\si^{j+1}))} \G^{{2\kappa}(\frac{5}{3}\kappa^{-1})^{j+1}}dyds \right)^{\frac{1}{2\kappa}(\frac{3}{5}\kappa)^{j+1}} \nn\\
&\lesssim&\frac{(1+E_R(b))^{\frac{10}{3\kappa}\sum^j_{i=0}(\frac{3}{5}\kappa)^{i+1}}}
{\left[\si^{\sum^j_{i=0}(i+1)(\frac{3}{5}\kappa)^{i+1}}R^{\sum^j_{i=0}(\frac{3}{5}\kappa)^{i+1}}\right]^{\frac{5}{3}(12-\frac{20}{p})\frac{1}{2\kappa}}}\times  \nn \\
&&\left(\int\int_{P(R)}\G^{2\kappa}dyds\right)^{\frac{1}{2\kappa}}.  \nn
\eea
Note that $\frac{3}{5}\kappa\in [\frac{3}{5},1)$ when we assume $p \in (\frac{5}{3},2]$. So all the sums on the above are convergent, let $j\rightarrow \infty$ yield that
\bea
\sup\limits_{(x,t)\in P(\frac{R}{2})}|\G|&\lesssim&\frac{\left(1+E_R(b)\right)^{\frac{15p-20}{6p-10}}}{R^{\frac{5}{2\kappa}}}\left(\int\int_{P(R)}\G^{{2\kappa}}dyds \right)^{\frac{1}{2\kappa}}  \nn  \\
&\lesssim&\left(1+E_R(b)\right)^{\frac{15p-20}{6p-10}}\left(\int\int_{P(R)}\frac{1}{R^5}\G^{{2\kappa}}dyds \right)^{\frac{1}{2\kappa}}  \nn  \\
&\lesssim&\left(1+E_R(b)\right)^{\frac{15p-20}{6p-10}}\left(\int\int_{P(R)}\frac{1}{R^5}\G^{\frac{2p}{3p-4}}dyds\right)^{\frac{3p-4}{2p}}.\label{e2.11}
\eea

Next we use \eqref{e2.11} and an algebraic trick to improve our estimate \eqref{e2.11}. This is from Li-Schoen \cite{Li01}.
From the process of proving \eqref{e2.11}, we have for $\g\in(0,\frac{1}{2}],\ \th\in[\frac{1}{2},1-\g]$
\bea
&&\sup\limits_{P(\th R)}\G^{\frac{2p}{3p-4}} \nn  \\
&\lesssim& (1+E_R(b))^{\frac{5p}{3p-5}}\frac{1}{R^5}\int\int_{P((\th+\g)R)}\G^{\frac{2p}{3p-4}}  \nn  \\
&\lesssim& (1+E_R(b))^{\frac{5p}{3p-5}}\frac{1}{R^5}\sup\limits_{P((\th+\g) R)}\G^{\frac{2p}{3p-4}-2}\int\int_{P((\th+\g) R)} \G^2.  \nn
\eea
Let $K\triangleq \frac{1}{R^5}\int\int_{P(R)}\G^2$,
then we have
\[
\sup\limits_{P(\th R)}\G^{\frac{2p}{3p-4}}\lesssim\left(1+E_R(b)\right)^{\frac{5p}{3p-5}}K\left(\sup\limits_{P((\th+\g) R)}\G^{\frac{2p}{3p-4}}\right)^{1-\frac{3p-4}{p}}.
\]
Define $M(\th)=\sup\limits_{P(\th R)}\G^{\frac{2p}{3p-4}}$,then we yield that
\[
M(\th)\lesssim(1+E_R(b))^{\frac{5p}{3p-5}}KM(\th+\g)K^\lambda,
\]
where $\lambda=1-\frac{3p-4}{p}$.
Choosing $\th_0=\frac{1}{2},\ \th_i=\th_{i-1}+\frac{1}{2^{i+1}}\ {\rm and}\  \g=\frac{1}{2^{i+1}}$,
then we get
\[
M(\th_0)\lesssim \left[K(1+E_R(b))^{\frac{5p}{3p-5}}\right]^{\sum^j_{i=1}\lambda^{i-1}}M(\th_j)^{\lambda^j}.
\]
For $\lambda<1$, letting $j\rightarrow \infty$, then we have
\[
M(\th_0)\lesssim \left[K(1+E_R(b))^{\frac{5p}{3p-5}}\right]^{\frac{p}{3p-4}}.
\]
That is
\[
\sup\limits_{P(\frac{1}{2}R)} \G^{\frac{2p}{3p-4}} \lesssim \left[K(1+E_R(b))^{\frac{5p}{3p-5}}\right]^{\frac{p}{3p-4}}.
\]
So
\bea
&&\sup\limits_{P(\frac{1}{2}R)} |\G|   \nn  \\
&\lesssim& (1+E_R(b))^{\frac{5p}{2(3p-5)}}K^{\frac{1}{2}}  \nn  \\
&\lesssim& \left(1+E_R(b)\right)^{\frac{5p}{2(3p-5)}}\left(\int\int_{P(R)}\frac{1}{R^5}\G^2dyds\right)^{\frac{1}{2}}. \nn
\eea
This proves our Lemma.  \ef
\section{Proof of Theorem 1.2}

In this section we study the continuity of $\G$ using the local maximum estimates \eqref{e2.9} and Nash type method for parabolic equations. First let us introduce some notations.

For $0< R\leq 1$, we define
\[
m_R=\inf\limits_{P(R)} \G, \q   M_R=\sup\limits_{P(R)}\G,\q J_R=M_R-m_R.
\]
Define
\begin{equation}
u=\left\{
\begin{aligned}
\frac{2(M_R-\G)}{J_R}             &\q    \textrm{if}\ M_R>-m_R,\\
\frac{2(\G-m_R)}{J_R}             &\q    \textrm{else} , \\
\end{aligned}
\right.  \label{e3.1}
\end{equation}
hence
\begin{equation}
0\leq u\leq 2,\q a\triangleq u|_{r=0}\geq 1. \label{e3.2}
\end{equation}
\\

\noindent{\bf Lower bound on $\|u\|_{L^q}$.}

We give a lemma to state that there is a lower bound on $\|u\|_{L^q}$ where $q\in(0,1)$. This bound depends on our $E_R(b)$ norm and will serve as an input for Nash's argument as we will describe it later on.

\begin{lemma}\label{l3.1}
If $u$ is a solution of \eqref{e1.30} and satisfies \eqref{e3.2}. Then for $\forall q\in (0,1)$ ,we have
\begin{equation}
\frac{1}{R^{\frac{5}{q}}}\|u\|_{L^q(P(\frac{R}{2}))}\gtrsim a(1+E_R(b))^{-\frac{8}{q}}.\label{e3.3}
\end{equation}
\end{lemma}
\noindent{\bf Proof.} Since the lemma is scaling invariant, we just take $R=1$ in the proof. Let $\psi(x,t)=\phi(x)\eta(t)$, where $\phi\in C^\infty_0$ s.t.
$\phi=1$ in $B_{\frac{1}{2}}$ and $\phi=0$ in $B^c_{1}$. $\frac{\nabla\phi}{\sqrt{\phi}}$ and $\nabla(\frac{\nabla\phi}{\sqrt{\phi}})$ are bounded. $\eta\in C^\infty_0$ s.t. $\eta=1$ in $[-\frac{7}{8},-\frac{1}{8}]$ and $\eta$ is supported in $(-1,0)$. Let us test \eqref{e1.30} by $qu^{q-1}\psi^2$, where
$q\in (0,\frac{1}{2})$.Then we have
\be
\int\int(\p_s u^q+b\cdot\nabla u^q+\frac{2}{r}\p_ru^q)\psi^2dyds=q\int\int\Delta u u^{q-1}\psi^2dyds. \label{e3.4}
\ee
Similarly as in \cite{Lei01}, we have
\bea
-\int\int\frac{2}{r}\p_ru^q\psi^2dyds&=&-\int\int2\p_ru^q\psi^2drdzds    \nn  \\
                                     &=&\int^0_{-1}\int2u^q\psi^2|_{r=0}dzds+\int\int \frac{4}{r}u^q \psi\p_r\psi dyds  \nn \\
                                     &\geq&-C\int\int u^qdyds+\frac{3}{2}a^q.  \label{e3.5}
\eea
Here we note that $\frac{\p_r \psi}{r}=\frac{\p_\rho \psi}{\rho}$ for $\phi$ is a radial function. Because $\phi=1$ near $\rho=0$, $\frac{\p_r \psi}{r}$ has no singularity.

Moreover
\bea
&&\int\int(-\p_s u^q+q\Delta uu^{q-1})\psi^2 dyds   \nn   \\
&=&\int\int2u^q[\psi\p_s\psi+|\nabla\psi|^2-\frac{q-2}{q}\psi\Delta\psi]dyds-\frac{4(q-1)}{q}\int\int |\nabla(u^{\frac{q}{2}}\psi)|^2 dyds   \nn  \\
&\geq&-C\int\int u^q dyds-\frac{4(q-1)}{q}\int\int|\nabla(u^{\frac{q}{2}}\psi)|^2 dyds.        \label{e3.6}
\eea
For the term involving $b$,we compute the same as \eqref{e2.7}
\bea
&&-\int\int b\cdot\nabla u^q\psi^2 dyds   \nn   \\
&\geq& -CE^4_R(b)\left(\int\int u^{\frac{q}{2}\frac{2p}{3p-4}}dyds\right)^{\frac{3p-4}{p}}-\frac{4(q-1)}{q}\int\int|\nabla(u^{\frac{q}{2}}\psi)|^2.  \label{e3.7}
\eea
Combining \eqref{e3.4},\eqref{e3.5},\eqref{e3.6} and \eqref{e3.7},we derive
\[
\int\int u^q dyds +E^4_R(b)\left(\int\int u^{\frac{q}{2}\frac{2p}{3p-4}}dyds\right)^{\frac{3p-4}{p}}\gtrsim a^q.
\]
Using H\"{o}lder inequality,we have
\[
\left(\int\int u^{2q}dyds\right)^{\frac{1}{2}}+E^4_R(b)\left(\int\int u^{2q}dyds\right)^{\frac{1}{2}}\gtrsim a^q,
\]
then we have
\[
(1+E^4_R(b))^{\frac{1}{q}}\left(\int\int u^{2q}\right)^{\frac{1}{2q}}\gtrsim a.
\]
So
\[
\left(\int\int u^{2q}dyds\right)^{\frac{1}{2q}}\gtrsim a(1+E_R(b))^{-\frac{4}{q}}.
\]
This proves our lemma, since $q\in(0,\frac{1}{2})$ is arbitrary.
\\

\noindent{\bf Nash's lower bound}

Before proving the Nash's lower bound estimates, we recall a Nash inequality ,whose proof can be found in \cite{Chen02}.
\begin{lemma}\label{l3.2}
Let $M\geq1$ be a constant and $\mu$ be a probability measure. Then for all $0\leq f\leq M$, there holds
\[
|\ln\int fd\mu -\int \ln f d\mu|\leq\frac{M\|g\|_{L^2}}{\int f d\mu}
\]
where $g=\ln f-\int \ln f d\mu$.
\end{lemma}

Now we come to prove Nash's lower bound estimate. We define a Lipschitz continuous cut-off function such that
\bea
&\zeta=1\  \rm{in}\ B(\frac{1}{2}),\ \zeta=0\ \rm{in}\  B(1)^c,\  \int_{\bR^3}\zeta^2dx=1.    \nn
\eea
In fact we take
\be
\zeta=c\left \{
\begin{aligned}
& 1  & {\rm in}\ B(\frac{1}{2});\\
& 2(1-|x|)& \q {\rm in}\ B(1)/B(\frac{1}{2}).
\end{aligned}\right.\label{e3.8}
\ee
where $c$ is a constant to ensure $\int_{\bR^3}\zeta^2dx=1$. Let $\zeta_R(x)=\frac{1}{R^{\frac{3}{2}}}\zeta(\frac{x}{R})$.
\\
\begin{lemma}\label{l3.3}
Let $0\leq u\leq 2$ be a solution of \eqref{e1.30} in $P(R)$ which is assumed to satisfy
\be
\|u\|_{L^1{(P(\frac{R}{2}))}} \geq c_1(1+E_R(b))^{-8} R^5.  \label{e3.9}
\ee
Moreover, we assume that $u|_{r=0}$ is a constant bigger than 1, then there exists a $\tau >0$ such that
\be
-\int_{\bR^3}\ln u \zeta^2_R dx \lesssim (1+E_R(b))^{24}, \q \rm{for}\ -\tau R^2\leq t<0. \nn
\ee
\end{lemma}
\noindent{\bf Proof.} First, let us define $u_R(x,t)=u(Rx,R^2t),\ b_R(x,t)=Rb(Rx,R^2t)$. It is clear that $u_R(x,t)$ solves the equation
\be
\p_t u_R+b_R\cdot\nabla u_R+\frac{2}{r}\p_r u=\Delta u_R  \q \rm{in}\ \rm{P}(1)   \nn
\ee
and $0\leq u_R\leq 2$, $\|u_R\|_{L^1{(P(\frac{1}{2}))}} \geq c_1(1+E_1(b_R))^{-8}$. The estimate we are going to get is
\be
-\int_{\bR^3}\ln u_R \zeta^2 dx \lesssim (1+E_1(b_R))^{24}, \q \rm{for}\ -\tau \leq t<0.  \nn
\ee
For convenience,we shall drop all $R$ and the subscript from now on and set $R=1$. \\
Also denote
\be
E=E_1(b_R).  \nn
\ee
Let $v=-\ln u$. It is easy to see that $v$ solves the equation
\be
\p_s v+ b\cdot v+\frac{2}{r}\p_r v-\Delta v+|\nabla v|^2 =0. \label{e3.10}
\ee
Testing \eqref{e3.10} by $\zeta^2$, we have
\be
\int \p_s v\zeta^2 dx +\int |\nabla v|^2\zeta^2dx=\int[\Delta v-\frac{2}{r}\p_r v-b\cdot v]\zeta^2 dx. \label{e3.11}
\ee
Using the Cauchy-Schwartz inequality and integration by parts, we have
\bea
\int\Delta v\zeta^2 dx&=&-2\int \nabla v\cdot\nabla\zeta\zeta dx   \nn  \\
                       &\leq&\frac{1}{4}\int|\nabla v|^2\zeta^2+4\int|\nabla\zeta|^2 dx \nn  \\
                       &\leq&\frac{1}{4}\int|\nabla v|^2\zeta^2 +C.   \label{e3.12}
\eea
Let $\ov(t)=\int v(\cdot,t)\zeta^2dx$, by recalling the assumption that $u|_{r=0}$ is a non-zero constant and the weighted poincar\'{e} inequality
\be
\int|v-\ov|^2\zeta^2dx\leq C\int|\nabla v|^2\zeta^2dx, \label{e3.13}
\ee
one can estimate
\bea
-\int\frac{2}{r}\p_r v \zeta^2&=&-\int\frac{2}{r}\p_r(v-\ov)\zeta^2 dx  \nn  \\
                              &=&-4\pi\int^{+\infty}_{-\infty}\int^{+\infty}_0\p_r(v-\ov)\zeta^2drdz  \nn  \\
                              &=&-4\pi\int^{+\infty}_{-\infty}\zeta^2(v-\ov)|^{\infty}_{r=0}dz+4\pi\int^{+\infty}_{-\infty}(v-\ov)2\zeta\p_r \zeta drdz  \nn \\
                              &=&4\pi\int^{+\infty}_{-\infty}v\zeta^2|_{r=0}dz-4\pi\ov\int^{+\infty}_{-\infty}\zeta^2dz+4\int(v-\ov)\zeta\frac{\p_r\zeta }{r}rdrd\th dz                                      \nn  \\
                              &\leq&C-C\ov(s)+\frac{1}{4}\int|\nabla v|^2\zeta^2 dx.     \label{e3.14}
\eea

Before estimating the term involving $b$, we need a more general weighted poincar\'{e} inequality.

Let $B_R$ be a ball centered at $0$ in $\bR^n$. Let $1\leq r\leq q<\i$ satisfy $\f{1}{q}\geq \f{1}{r}-\f{1}{n}$ and $1-\f{n+2}{r}+\f{n+2}{q}\leq 0$, then we have
\be
\left(\frac{1}{|B_R|}\int_{B_R}(v-\ov)^{q}\zeta^2_R dx\right)^{\frac{1}{q}}\leq C R^{1-\f{2}{r}+\f{2}{q}} \left(\frac{1}{|B_R|}\int_{B_R}|\nabla v|^r\zeta^2_Rdx\right)^{\frac{1}{r}}, \label{e3.15}
\ee
here $|B_R|$ means the Lebesgue measure of the ball $B_R$ and $C$ depends only on $q$, $r$, $n$. One can see \cite{Chua} for its proof.
Hence, due to the divergence-free of $b$ and H\"{o}lder inequality, we have
\bea
-\int(b\cdot\nabla)v \zeta^2 dx&=&\int 2\zeta(v-\ov)b\cdot\nabla\zeta dx  \nn   \\
                               &\leq&\left(\int|b|^p|\nabla\zeta|^p dx\right)^{\frac{1}{p}}\left(\int(v-\ov)^{q}\zeta^{q} dx\right)^{\frac{1}{q}} \nn \\
                               &\leq&\left(\int|b|^p|\nabla\zeta|^p dx\right)^{\frac{1}{p}}\left(\int(v-\ov)^{q}\zeta^{2} dx\right)^{\frac{1}{q}},  \nn
\eea
here $\frac{1}{p}+\frac{1}{q}=1$.

In  \eqref{e3.15}, let $R=1$, $r=2,\ n=3$. When $p\in (\frac{5}{3}, 2]$, $q=\frac{p}{p-1}$ can satisfy $\f{1}{q}\geq \f{1}{r}-\f{1}{n}$ and $1-\f{n+2}{r}+\f{n+2}{q}\leq 0$. So
\be
\left(\int(v-\ov)^{q}\zeta^{2} dx\right)^{\frac{1}{q}}\leq C\left(\int|\nabla v|^2\zeta^2dx\right)^{\frac{1}{2}}. \nn
\ee
Hence
\bea
-\int(b\cdot\nabla)v \zeta^2 dx&\leq&C\left(\int|b|^p|\nabla\zeta|^p dx\right)^{\frac{1}{p}}\left(\int|\nabla v|^2\zeta^{2} dx\right)^{\frac{1}{2}} \nn \\
                               &\leq&C\left(\int_{B(1)}|b|^p dx\right)^{\frac{2}{p}} +\frac{1}{4}\int|\nabla v|^2\zeta^{2} dx \nn  \\
                               &\leq&CE^2 +\frac{1}{4}\int|\nabla v|^2\zeta^{2} dx.  \label{e3.16}
\eea
Combining \eqref{e3.11},\eqref{e3.12},\eqref{e3.14} and \eqref{e3.16}, we have
\[
\p_s\int v\zeta^2 dx+C\ov(s)\leq -\frac{1}{4}\int|\nabla v|^2\zeta^2dx+(1+E)^2.
\]

Now we apply the Nash inequality, taking $f=u,\  d\mu=\zeta^2dx$ in Lemma 3.2, one has
\[
|\ln\int u\zeta^2 dx+\int v\zeta^2dx|^2\left(\int u \zeta^2 dx\right)^2\leq M^2\int|-v+\int v\zeta^2dy|^2\zeta^2 dx,
\]
here $M=2$ is the upper bound of $u$. Using the weighted Poincar\'{e} inequality \eqref{e3.13} once again, we have
\[
|\ln\int u\zeta^2 dx+\int v\zeta^2dx|^2\left(\int u \zeta^2 dx\right)^2\leq C\int |\nabla v|^2\zeta^2 dx.
\]
Then finally we obtain
\[
\p_s\ov(s)+C_0\ov(s)\leq(1+E)^2-\frac{1}{4C}|\ln\int u\zeta^2dx+\ov(s)|^2\left(\int u \zeta^2dx\right)^2.
\]
Recalling \eqref{e3.9}
\[
\|u\|_{L^1(P(\frac{1}{2}))}\geq c_1(1+E)^{-8}.
\]
Let $\chi$ be the characteristic function of the non-empty set
\[
W=\left\{s\in [-\frac{1}{4},0]:\|u\|_{L^1(B(\frac{1}{2}))}\geq \frac{c_1(1+E)^{-8}}{10}\right\}.
\]
We assert that $|W|\geq\frac{c_1(1+E)^{-8}}{20}$. In fact, if $|W|\leq\frac{c_1(1+E)^{-8}}{20}$, then
\bea
\|u\|_{L^1(P(\frac{1}{2}))}&<&\int_W 2|B(\frac{1}{2})|ds+\int_{W^c}\frac{c_1(1+E)^{-8}}{10}ds \nn  \\
                           &<&\frac{\pi}{3}|W|+\frac{c_1(1+E)^{-8}}{40}        \nn        \\
                           &<&(\frac{\pi}{60}+\frac{1}{40})c_1(1+E)^{-8}        \nn      \\
                           &<&c_1(1+E)^{-8},         \nn
\eea
this is a contradiction with \eqref{e3.9}.
Thus,one has
\bea
\p_s\ov(s)+C_0\ov(t)&\leq&(1+E)^2-\frac{1}{4C}\chi(s)|\ln\int u\zeta^2dx+\ov(s)|^2\left(\int u\zeta^2 dx\right)^2  \nn  \\
                  &\leq&(1+E)^2-\frac{1}{4C}\chi(s)|\ln\int u\zeta^2dx+\ov(s)|^2 \frac{c^2_1(1+E)^{-16}}{100}.  \label{e3.17}
\eea
The last inequality is due to
\be
\left(\int u\zeta^2 dx\right)^2\geq \left(\int_{B(\frac{1}{2})} u dx\right)^2\geq\frac{c^2_1(1+E)^{-16}}{100},     \nn
\ee
when $s\in W$.\\
From\eqref{e3.17}, we first have $\p_s\ov(s)+C_0\ov(s)\leq(1+E)^2$. This gives, for $-\frac{1}{4}\leq s_1\leq s_2\leq 0$,
\bea
\ov(s_2)&\leq& e^{C_0(s_1-s_2)}\ov(s_1)+e^{-s_2}\int^{s_2}_{s_1}(1+E)^2ds  \nn  \\
        &\leq& e^{C_0|s_1-s_2|}\ov(s_1)+C(1+E)^2.   \label{e3.18}
\eea
Now we consider two cases.\\
Case one: if there exists some $s_0\in [-\frac{1}{4},-\frac{c_1}{40}]$, such that
\be
\ov(s_0)\leq \frac{2}{C_0}(1+E)^2+4|\ln\frac{10}{c_1(1+E)^{-8}}|.   \nn  \\
\ee
Then for $s\in(s_0,0)$, from \eqref{e3.18},we have
\bea
\ov(s)&\lesssim& \ov(s_0)+(1+E)^2  \nn   \\
      &\lesssim& (1+E)^2+\ln \frac{10}{c_1}+\ln(1+E)  \nn  \\
      &\lesssim& (1+E)^2.  \nn
\eea
Choosing $\tau=s_0$, this completes the proof of the Lemma.\\
Case two:
if for any $s\in[-\frac{1}{4},-\frac{c_1}{40}]$,
\be
\ov(s)\geq \frac{2}{C_0}(1+E)^2+4|\ln\frac{10}{c_1(1+E)^{-8}}|.  \nn
\ee
Then when $s\in W\cap[-\frac{1}{4},-\frac{c_1}{40}]$,
\be
\ln \int u\zeta^2dx \geq \ln\int_{B(\frac{1}{2})}u\zeta^2\geq \ln\frac{c_1(1+E)^{-8}}{10}.  \nn
\ee
So we have
\be
\ov(s)+\ln\int u\zeta^2dx \geq C \ov(s).           \nn
\ee
From \eqref{e3.17},we have
\bea
\p_s \ov(s)+C_0\ov(s)\lesssim-C\chi(s)\frac{c^2_1(1+E)^{-16}}{100}\ov^2(s).   \nn
\eea
Then integrating the above inequality from $[-\frac{1}{4},-\frac{c_1}{40}]$, one  gets
\bea
\ov(-\frac{1}{4})^{-1}-\ov(-\frac{c_1}{40})^{-1}&\lesssim& -c^2_1(1+E)^{-16}\int \chi ds   \nn  \\
                                                &\lesssim& -c^2_1(1+E)^{-16}|W|    \nn   \\
                                                &\lesssim& -c^3_1(1+E)^{-24}.  \nn
\eea
Since $\ov(-\frac{1}{4})\geq 0$ in this case, we have
\be
\ov({-\frac{c_1}{40}})\lesssim \frac{(1+E)^{24}}{c^3_1}.  \nn
\ee
Then we use \eqref{e3.18},for $s\in[-\frac{c_1}{40},0]$,
\be
\ov(s)\lesssim(1+E)^2+(1+E)^{24}\lesssim(1+E)^{24}.   \nn
\ee
So we can take $\tau=\frac{c_1}{40}$, this proves the lemma.  \ef
As a corollary of Lemma 3.3, we derive a lower bound of positive solution of \eqref{e1.30}.
\begin{corollary}
Let $u,\ \tau$ be given in Lemma 3.3 and $E_R(b)$ satisfies the assumption \eqref{e1.6}. Then there exists a $\delta\in (0,1)$, depending only on $R$,
such that
\be
\inf\limits_{P(\frac{\sqrt{\tau}}{2}R)} u\geq \frac{1}{2}\delta(R).  \label{e3.19}
\ee
In fact, we take  $\dl=|\ln\f{R}{3}|^{-1}$.
\end{corollary}
\noindent{\bf Proof.}Using Lemma\eqref{l3.3}, one has
\bea
(1+E_R(b))^{24}&\gtrsim&-\int\zeta^2_R(x)\ln u(x,t)dx   \nn  \\
              &=&-\int_{\delta<u\leq 1}\zeta^2_R(x)\ln udx-\int_{u\leq \delta}\zeta^2_R(x)\ln udx   \nn  \\
              &&\int_{1<u\leq 2}\zeta^2_R(x)\ln udx     \nn   \\
              &\geq&-\int_{u\leq \delta}\zeta^2_R(x)\ln udx-\ln2\int_{1<u\leq 2}\zeta^2_R(x)dx  \nn   \\
              &\geq&-\int_{u\leq \delta}\zeta^2_R(x)\ln udx-\ln2 .                          \nn
\eea
This implies that
\be
-\int_{u\leq \delta}\zeta^2_R(x)\ln udx\lesssim(1+E_R(b))^{24}.  \nn
\ee
For $t\in[-\tau R^2,0]$, consequently, one has
\be
|\left\{x\in B(\frac{R}{2})|u\leq \delta\right\}|\lesssim \frac{R^3}{-\ln \delta}(1+E_R(b))^{24}.  \nn
\ee
Using the mean value inequality \eqref{e2.9}, one has
\bea
\sup\limits_{P(\frac{\s{\tau}}{2}R)}(\delta-u)_+&\lesssim&(1+E_R(b))^{\frac{5p}{2(3p-5)}}\left(\int\int_{P(\s{\tau} R)}\frac{1}{(\s{\tau} R)^5}(\delta-u)^2_+ dy ds\right)^{\frac
                                                     {1}{2}}  \nn   \\
                                            &\lesssim&(1+E_R(b))^{\frac{5p}{2(3p-5)}}\delta\left[\frac{(\s{\tau} R)^2}{(\s{\tau} R)^5}\frac{R^3}{-\ln \delta}(1+E_R(b))^{24}\right]^{\frac{1}{2}}                \nn  \\
                                            &\lesssim&\left(1+E_R(b)\right)^{\frac{5p}{2(3p-5)}+12}\frac{\delta}{\sqrt{-\ln\delta}}.  \nn
\eea
This gives
\be
\inf\limits_{P(\frac{\s{\tau}}{2}R)} u\geq \delta\left[1-C\frac{(1+E_R(b))^{\frac{5p}{2(3p-5)}+12}}{\sqrt{-\ln\delta}}\right].  \nn
\ee
Under the assumption \eqref{e1.6}, one has
\bea
\inf\limits_{P(\frac{\s{\tau}}{2}R)} u&\geq& \delta\left\{1-C\frac{\left[\left(\ln\ln\frac{3}{R}\right)^{\frac{3p-5}{77p-120}(1-\b)}\right]^{\frac{77p-120}{6p-10}}}{\sqrt{-\ln\delta}}\right\} \nn \\
                                  &\geq& \delta\left[1-C\frac{\left(\ln\ln\frac{3}{R}\right)^{\frac{1-\b}{2}}}{\sqrt{-\ln\delta}}\right],  \nn
\eea
when $R\in (0,1)$. We can take $\delta(R)=(\ln\frac{3}{R})^{-1}$ to ensure $\inf\limits_{P(\frac{\s{\tau}}{2}R)} u\geq\frac{1}{2}\delta(R)$. This proves the corollary.\ef
\noindent{\bf Proof of theorem 1.2}

We define
\be
m_\tau= \inf\limits_{P(\frac{\s{\tau}}{2}R)} \G, \q M_\tau=\sup\limits_{P(\frac{\s{\tau}}{2}R)} \G.  \nn
\ee
Then from \eqref{e3.1} and \eqref{e3.19}, one has
\be
\f{1}{2}\dl(R)\leq\inf\limits_{P(\frac{\s{\tau}}{2}R)}u=\left \{
\begin{aligned}
&2(M_R-M_\tau)/J_R   \q\rm{if}\  M_R>-m_R;\\
&2(m_\tau-m_R)/J_R   \q\rm{else}.    \\
\end{aligned}
\right.\nn
\ee
We add the two cases together to get that
\be
\delta(R )\leq\frac{4}{J_R}\left\{J_R-\rm{osc}(\G,\frac{\s{\tau}}{2}R)\right\},  \nn
\ee
here $\rm{osc}(\G,\frac{\s{\tau}}{2}R)=M_\tau-m_\tau$ and $\rm{osc}(\G,R)=M_R-m_R=J_R$. So
\be
\rm{osc}(\G,\frac{\s{\tau}}{2}R)\leq\left(1-\frac{\delta(R)}{4}\right)\rm{osc}\left(\G,R\right).  \nn
\ee
We write it also as
\be
J_{\frac{\s{\tau}}{2}R}\leq \left(1-\frac{\delta(R)}{4}\right)J_{R}.  \nn
\ee
For any small $R$, there exists some integer $j\geq 0$ such that $(\frac{\s{\tau}}{2})^{j+1}<R/3\leq(\frac{\s{\tau}}{2})^j$. Using the above inequality, an
iteration argument gives
\be
J_{R/3}\leq J_{(\frac{\s{\tau}}{2})^j}\leq\prod^{j-1}_{k=0}\left[1-\frac{\delta((\frac{\s{\tau}}{2})^{k})}{4}\right]J_1. \nn
\ee
Noting that $\ln(1-x)\leq -x$ for sufficiently small positive $x$, one has
\bea
J_{R/3}&\leq&\exp\ln J_{R/3}                                   \nn    \\
   &\leq&J_1\exp\sum^{j-1}_0\ln\left[1-\frac{\delta\left((\frac{\s\tau}{2})^k\right)}{4}\right]  \nn \\
   &\leq&J_1\exp\left\{-\frac{1}{4}\sum^{j-1}_0\delta((\frac{\s\tau}{2})^k)\right\}     \nn  \\
   &\leq&J_1\exp\left\{-\frac{1}{4}\sum^{j-1}_0\left[\ln\frac{3}{(\frac{\s{\tau}}{2})^k}\right]^{-1}\right\}  \nn  \\
   &\leq&J_1\exp\left\{-\frac{1}{4}\sum^{j-1}_0\left[\ln3+k|\ln(\frac{\s{\tau}}{2})|\right]^{-1}\right\}  \nn  \\
   &\leq&J_1\exp\left\{-\frac{1}{4}\sum^{j-1}_0\left[C(k+1)\right]^{-1}\right\}     \nn   \\
   &\lesssim&J_1\exp\left\{-c_0\ln j\right\}      \nn    \\
   &\lesssim&J_1j^{-c_0}            \nn    \\
   &\lesssim&J_1|\ln \f{R}{3}|^{-c_0}.   \nn
\eea
Since, $\G|_{r=0}=0$, the above estimate proves our Theorem 1.2. \ef

\section{Proof of Theorem 1.1}

In this section we will prove Theorem 1.1 and get the regularity of the solution under the assumption \eqref{e1.2}. The idea comes from [Chen-Strain-Tsai-Yau]'s proof where they assume $|v|\leq Cr^{-1}$.\\
We divide the proof into 3 steps.\\
{\bf Step one: scaling of the solution and set up of a equation}

Let $M$ be the maxmium of $|v|$ up to a fixed time $t_0$ and we may assume $M>1$ is large. Define the scaled solution
\be
v^M(X,T)=M^{-1}v(\frac{X}{M},\frac{T}{M^2}),\q X=(X_1,X_2,Z).  \nn
\ee
Denote $x=(x_1,x_2,z)$ and $X=(X_1,X_2,Z)$, $r=\s{x^2_1+x^2_2}$ and $R=\s{X^2_1+X^2_1}$. We have the following estimate for $r$ and $R$ for time $t<t_0$ and $T< M^2t_0$:
\be
|\nabla^k v^M|\leq C_k.   \label{e4.1}
\ee
This inequality follows from $\|v^M\|_{L^\infty}\leq1$ for $t<t_0$ and the standard regularity theorem of Navier-Stokes equations. Its angular component(we omit
the time dependence below) $v^M_\th(R,Z)$ satisfies $v^M_\th(0,Z)=0=\p_Z v^M_\th(0,Z)$ for all $Z$. By mean value theorem and \eqref{e4.1},
\be
|v^M_\th(R,Z)|\lesssim R,\ |\p_Z v^M_\th(R,Z)|\lesssim R \q \rm{for} \ R\leq1. \nn
\ee
Together with \eqref{e4.1} for $R\geq1$, we get
\be
|v^M_\th|\lesssim \min\{R,1\}, \  |\p_z v^M_\th|\lesssim \min\{R,1\}.  \label{e4.2}
\ee
Due to Theorem 1.2,
\be
|\G|=|rv_\th(r,z)|\lesssim
\left\{
\begin{aligned}
(\ln\frac{3}{r})^{-c_0} \ &\q {\rm for}\ r\leq1;  \\
1\q                       &\q {\rm for}\ r>1.
\end{aligned}
\right.\nn
\ee
That is
\be
|v_\th(r,z)|\lesssim
\left\{
\begin{aligned}
\frac{(\ln\frac{3}{r})^{-c_0}}{r} \ &\q{\rm for}\ r\leq1;  \\
\frac{1}{r}\q                       &\q {\rm for}\ r>1.
\end{aligned}
\right.\label{e4.3}
\ee
Then $v^M_\th(R,Z)$ satisfies the estimate
\be
|v^M_\th(R,Z)|=M^{-1}|v_\th(\frac{X}{M},\frac{T}{M^2})|\lesssim
\left\{
\begin{aligned}
\frac{(\ln\frac{3M}{R})^{-c_0}}{R} \ &\q {\rm for}\ R\leq M;  \\
\frac{1}{R}\q                       &\q {\rm for}\ R>M.
\end{aligned}
\right.\nn
\ee
Combining this with \eqref{e4.2}, one has
\be
|v^M_\th(R,Z)|\lesssim
\left\{
\begin{aligned}
\min\left\{R,\frac{(\ln\frac{3M}{R})^{-c_0}}{R}\right\} \ &\q {\rm for}\ R\leq 1;  \\
\min\left\{1,\frac{(\ln\frac{3M}{R})^{-c_0}}{R}\right\} \     &\q {\rm for}\ 1<R<M; \\
\frac{1}{R}\q\q&\q {\rm for}\ R\geq M.
\end{aligned}
\right.\label{e4.4}
\ee

Now consider the angular component of the rescaled vorticity. Recall $\Omega=\frac{w_\th}{r}$. Let
\bea
f&=&\Omega^M(X,T)=M^{-3}\Omega(\frac{X}{M},\frac{T}{M^2})=M^{-3}w_\th(\frac{X}{M},\frac{T}{M^2})\frac{M}{R}  \nn  \\
 &=&\frac{w^M_\th(X,T)}{R},   \nn
\eea
where
\be
w^M_\th(X,T)=w_\th(\frac{X}{M},\frac{T}{M^2})M^{-2}.    \nn
\ee
Note that $w^M_\th$ and $\nabla w^M_\th$ are bounded by \eqref{e4.1} and also $w^M_\th|_{R=0}=0$, so one has
\be
|f|\lesssim\frac{1}{1+R}.  \nn
\ee
From the equation \eqref{e1.4}, $f$ satisfies
\be
(\p_T-L)f=g,\q L=\Delta+\frac{2}{R}\p_R-b^M\cdot\nabla,     \nn
\ee
where $g=R^{-2}\p_Z(v^M_\th)^2$ and $b^M=v^M_Re_R+v^M_Ze_Z$, $|b^M|\leq 1$.\\
Combining the estimates \eqref{e4.2} and \eqref{e4.4}, one has
\be
g=\frac{2}{R^2}v^M_\th\p_Z(v^M_\th)\lesssim
\left\{
\begin{aligned}
\min\left\{1,\frac{(\ln\frac{3M}{R})^{-c_0}}{R^2}\right\}\ \q &\q {\rm for}\ R\leq 1;  \\
\min\left\{\frac{1}{R^2},\frac{(\ln\frac{3M}{R})^{-c_0}}{R^3}\right\} \q    &\q {\rm for}\ 1<R<M; \\
\frac{1}{R^3}\q\q&\q {\rm for}\ R\geq M.
\end{aligned}
\right.\label{e4.5}
\ee
Let $P(X,T;Y,S)$ be the kernel of $\p_T-L$. By Duhamel's formula
\bea
f(X,T)&=&\int P(X,T;Y,S)f(Y,S)dY+\int^T_S\int P(X,T;Y,\tau)g(Y,\tau)dYd\tau  \label{e4.60} \nn  \\
      &:=&I+II.  \label {4.6}
\eea
{\bf Step two: bounding of $f$ }
\\
In the following ,we will estimate \eqref{e4.60} and give a bound for $f(X,T)$.

The kernel $P(X,T;Y,S)$ satisfies $P\geq0,\ \int P(X,T;Y,S) dY\leq1$ and
\be
P(X,T;Y,S)\leq C(T-S)^{-3/2}\exp\left\{-C\frac{|X-Y|^2}{T-S}\left(1-\frac{T-S}{|X-Y|}\right)^2_+\right\}.\label{e4.6}
\ee

The proof of estimate \eqref{e4.6} is based on \cite{Carlen01}, but due to the singularity of the term $\f{2}{r}\p_r$, the proof is more involved. For completeness of our paper, we will prove it in Section 5 as Theorem1.3.\\
Now we give estimates of $P$ in two cases. \\
Case one: when  $1-\frac{T-S}{|X-Y|}>\frac{1}{2}$, that is $|X-Y|>2(T-S)$,
\bea
\exp\left\{-\frac{|X-Y|^2}{T-S}\left(1-\frac{T-S}{|X-Y|}\right)^2_+\right\}&\leq& \exp\left\{-\frac{1}{4}\frac{|X-Y|^2}{T-S}\right\}  \nn  \\
                                                 &\leq&\left\{
                                                 \begin{aligned}
                                                 \exp\left\{-\frac{1}{4}\frac{|X-Y|}{T-S}\right\}    \q&\q {\rm for}\ |X-Y|\geq1;\\
                                                 \exp \left\{-\frac{1}{4}\frac{|X-Y|^2}{T-S}\right\}  \q&\q {\rm for}\ |X-Y|<1.
                                                 \end{aligned}
                                                 \right.\nn
\eea
Case two: when $1-\frac{T-S}{|X-Y|}\leq\frac{1}{2}$, that is $|X-Y|\leq2(T-S)$,
\be
\exp\left\{-\frac{|X-Y|^2}{T-S}\left(1-\frac{T-S}{|X-Y|}\right)^2_+\right\}\leq1\leq e^2\exp\left\{-\frac{|X-Y|}{T-S}\right\}. \nn
\ee
With these estimates and H\"{o}lder inequality, one gets, for $I$ in \eqref{e4.60},
\bea
|I|&\leq&\left[\int P(X,T;Y,S)|f(Y,S)|^3dY\right]^{\frac{1}{3}}  \left[\int P(X,T;Y,S)dY\right]^{\frac{2}{3}}\nn \\
   &\leq&\left\{\int^{+\infty}_{-\infty}\int^{+\infty}_0(T-S)^{-\frac{3}{2}}\left[e^{-\frac{|X_3-Y_3|}{T-S}}+e^{-\frac{|X_3-Y_3|^2}{T-S}}\right]
   \frac{R}{(R+1)^3}dRdY_3\right\}^{1/3} \nn  \\
   &\lesssim&(T-S)^{-\frac{1}{2}}\left\{\int^{+\infty}_{-\infty}\left[e^{-\frac{|X_3-Y_3|}{T-S}}+e^{-\frac{|X_3-Y_3|^2}{T-S}}\right]dY_3\right\}^{1/3} \nn \\
   &\lesssim&(T-S)^{-\frac{1}{2}}\left\{(T-S)+(T-S)^{\frac{1}{2}}\right\}^{1/3}    \nn \\
   &\lesssim&(T-S)^{-\frac{1}{6}}  \label{e4.7}
\eea
for $T-S\geq 1$, next
\bea
|II|&\leq&\int^T_S(T-\tau)^{-\frac{3}{2}}\left\{\int_{|X-Y|\leq2(T-\tau)}e^{-\frac{|X-Y|}{T-\tau}}|g|dY\right.   \nn  \\
    &&\left.+\int_{|X-Y|\geq2(T-\tau),|X-Y|>1}e^{-\frac{1}{4}\frac{|X-Y|}{T-\tau}}|g|dY
    +\int_{|X-Y|\geq2(T-\tau),|X-Y|<1}e^{-\frac{1}{4}\frac{|X-Y|^2}{T-\tau}}|g|dY\right\}d\tau  \nn  \\
    &:=&\int^T_S(T-\tau)^{-\frac{3}{2}}\left\{II_1+II_2+II_3\right\}d\tau.    \label{4.90}
\eea
We deal with $II_1,II_2,II_3$ in \eqref{4.90} as follows,
\bea
II_1+II_2&\lesssim&\int^{+\infty}_{-\infty}\int^{+\infty}_0\left(e^{-\frac{1}{4}\frac{|X_3-Y_3|}{T-\tau}}+e^{-\frac{|X_3-Y_3|}{T-\tau}}\right)|g|RdRdY_3 \nn \\
         &\lesssim&(T-\tau)\int^{+\infty}_0|g|RdR   \nn  \\
         &\lesssim&(T-\tau)\left\{\int^1_0\min\left\{R,\frac{(\ln\frac{3M}{R})^{-c_0}}{R}\right\}dR\right.    \nn \\
         &&\left.+\int^M_1\min\left\{\frac{1}{R},\frac{(\ln\frac{3M}{R})^{-c_0}}{R^2}\right\}dR +\int^{+\infty}_M\frac{1}{R^2}dR\right\}  \nn \\
         &:=&(T-\tau)(III_1+III_2+III_3) \label{4.10}
\eea
For $III_1$, when $R\in (0,1]$, the function $R$ is increasing while  $\frac{(\ln\frac{3M}{R})^{-c_0}}{R}$ is decreasing. Let $R_0$ be such that
\be
R_0=\frac{(\ln\frac{3M}{R_0})^{-c_0}}{R_0}.  \nn
\ee
This makes
\be
(\frac{1}{R_0})^{\frac{2}{c_0}}=\ln3M+\ln\frac{1}{R_0}\leq\ln3M+\frac{c_0}{2}(\frac{1}{R_0})^{\frac{2}{c_0}}. \nn
\ee
That is
\be
(1-\frac{c_0}{2})(\frac{1}{R_0})^{\frac{2}{c_0}}\leq\ln3M\leq(\frac{1}{R_0})^{\frac{2}{c_0}}. \nn
\ee
So, there exists a $C>1$, such that
\be
C^{-1}(\ln3M)^{-\frac{c_0}{2}}\leq R_0\leq C(\ln3M)^{-\frac{c_0}{2}}. \label{4.110}
\ee
Then
\be
\min\{R,\frac{(\ln\frac{3M}{R})^{-c_0}}{R}\}=
\left\{
\begin{aligned}
R                                 \q\q&\q{\rm for}\ 0\leq R\leq R_0; \\
\frac{(\ln\frac{3M}{R})^{-c_0}}{R}\q&\q{\rm for}\ R_0<R<1.
\end{aligned}
\right. \nn
\ee
By \eqref{4.110}, the control of $III_1$ in \eqref{4.10} is
\bea
III_1&\leq& \int^{R_0}_0 RdR+\int^1_{R_0}\frac{(\ln\frac{3M}{R})^{-c_0}}{R}dR  \nn  \\
     &\leq&\frac{1}{2}R^2_0+(\ln3M)^{-c_0}\int^1_{R_0}\frac{1}{R}dR  \nn  \\
     &\lesssim&(\ln3M)^{-c_0}(1+\ln\frac{1}{R_0})   \nn   \\
     &\lesssim&(\ln3M)^{-c_0}(1+\frac{1}{R_0})   \nn   \\
     &\lesssim&(\ln3M)^{-c_0}(1+(\ln3M)^{\frac{c_0}{2}})  \nn   \\
     &\lesssim&(\ln3M)^{-\frac{c_0}{2}}.   \label{e4.8}
\eea
For $III_2$ in \eqref{4.10}, one has
\bea
III_2&\leq&\int^M_1\frac{(\ln\frac{3M}{R})^{-c_0}}{R^2}dR  \nn  \\
     &=&M^{-1}\int^1_{M^{-1}}\frac{(\ln\frac{3}{R})^{-c_0}}{R^2}dR  \nn  \\
     &=&M^{-1}\left(\int^{M^{-\frac{1}{2}}}_{M^{-1}}+\int^1_{M^{-\frac{1}{2}}}\right)\frac{(\ln\frac{3}{R})^{-c_0}}{R^2}dR  \nn  \\
     &\leq&M^{-1}\left\{(\ln\frac{3}{M^{-\frac{1}{2}}})^{-c_0}\int^{M^{-\frac{1}{2}}}_{M^{-1}}\frac{1}{R^2}dR
     +\int^1_{M^{-\frac{1}{2}}}\frac{(R\ln\frac{3}{R})^{-c_0}}{R^{2-c_0}}dR\right\}.  \nn
\eea
Since $R\ln\frac{3}{R}$ is increasing when $R\in (0,1)$, one has,
\bea
III_2&\lesssim&M^{-1}\left\{(\ln3M)^{-c_0}M+(M^{-\frac{1}{2}}\ln3M)^{-c_0}M^{-\frac{c_0-1}{2}}\right\} \nn \\
     &\lesssim&(\ln3M)^{-c_0}+M^{-\frac{1}{2}}(\ln3M)^{-c_0}  \nn  \\
     &\lesssim&(\ln3M)^{-c_0}.\label{e4.9}
\eea
For $III_3$ in \eqref{4.10}, obviously
\be
III_3\lesssim M^{-1}. \label{e4.10}
\ee
Hence, Combining \eqref{e4.8}, \eqref{e4.9} and \eqref{e4.10}, from \eqref{4.10}, one has
\be
II_1+II_2\lesssim (T-\tau)(\ln3M)^{-\frac{c_0}{2}}.  \label{4.14}
\ee
For $II_3$ in \eqref{4.90}, using the Cauchy-Schartz inequality,one has
\bea
II_3&=&\int_{|X-Y|\geq2(T-\tau),|X-Y|<1}e^{-\frac{1}{4}\frac{|X-Y|^2}{T-\tau}}|g|dY   \nn  \\
    &=&\int e^{-\frac{1}{4}\frac{|X'-Y'|^2}{T-\tau}-\frac{1}{8}\frac{|X_3-Y_3|^2}{T-\tau}}e^{-\frac{1}{8}\frac{|X_3-Y_3|^2}{T-\tau}}|g|dY \nn  \\
    &\leq&\left(\int e^{-\frac{1}{2}\frac{|X'-Y'|^2}{T-\tau}-\frac{1}{4}\frac{|X_3-Y_3|^2}{T-\tau}}dY\right)^{\frac{1}{2}}
    \left(\int e^{-\frac{1}{4}\frac{|X_3-Y_3|^2}{T-\tau}}g^2dY\right)^{\frac{1}{2}}    \nn   \\
    &\leq&\left(\int e^{-\frac{1}{2}\frac{|X'-Y'|^2}{T-\tau}-\frac{1}{4}\frac{|X_3-Y_3|^2}{T-\tau}}dY\right)^{\frac{1}{2}}
    \left(\int^{+\i}_{-\i}\int^{+\i}_0 e^{-\frac{1}{4}\frac{|X_3-Y_3|^2}{T-\tau}}g^2RdRdY_3\right)^{\frac{1}{2}}    \nn   \\
    &\leq&(T-\tau)^{\frac{3}{4}}\cdot(T-\tau)^{\frac{1}{4}}\left(\int^{+\i}_0 \sup\limits_{Y_3}g^2RdR\right)^{\frac{1}{2}} \nn  \\
    &\leq&(T-\tau)\left(\int^{+\i}_0 \sup\limits_{Y_3}g^2RdR\right)^{\frac{1}{2}}. \nn
\eea
Following \eqref{e4.5} one has
\be
g^2\leq\left\{
\begin{aligned}
\min\left\{1,\frac{(\ln\frac{3M}{R})^{-c_0}}{R^2},\frac{(\ln\frac{3M}{R})^{-2c_0}}{R^4}\right\}\ \q &\q{\rm for}\ R\leq 1;  \\
\min\left\{\frac{1}{R^4},\frac{(\ln\frac{3M}{R})^{-c_0}}{R^5},\frac{(\ln\frac{3M}{R})^{-c_0}}{R^6}\right\} \q    &\q{\rm for}\ 1<R<M; \\
\frac{1}{R^6}\q\q\q&\q{\rm for}\ R\geq M. \nn
\end{aligned}
\right.
\ee
As the previous proof for \eqref{4.14}, one can get
\be
II_3\lesssim (T-\tau)(\ln3M)^{-\frac{c_0}{2}}. \label{4.15}
\ee
Inserting \eqref{4.14} and \eqref{4.15} into \eqref{4.90}, one has
\be
|II|\lesssim\int^T_S(T-\tau)^{-\frac{1}{2}}(\ln3M)^{-\frac{c_0}{2}}d\tau\lesssim(T-S)^{\frac{1}{2}}(\ln3M)^{-\frac{c_0}{2}}. \label{e4.11}
\ee
So, combining \eqref{e4.7} and \eqref{e4.11} ,from \eqref{4.6}, one has
\be
|f(X,T)| \lesssim(T-S)^{-\frac{1}{6}}+(T-S)^{\frac{1}{2}}(\ln3M)^{-\frac{c_0}{2}}.  \nn
\ee
Let $S=T-(\ln3M)^{-\frac{3}{4}c_0}>-M^2$(hence $f$ is defined),so
\be
|f(X,T)| \lesssim(\ln3M)^{-\frac{1}{8}c_0}.\nn
\ee
{\bf Step three: bounding the solution $v$ from $f$}\\
First
\be
|w_\th(x,t)|\leq M^2|w^M_\th(rM,zM,tM^2)|\leq |\Omega^M(rM,zM,tM^2)|M^2rM\leq CM^3r(\ln3M)^{-\frac{1}{8}c_0}. \nn
\ee
Therefore
\be
|w_\th(x,t)|\leq CM^2(\ln3M)^{-\frac{1}{16}c_0},\q {\rm for} \ r\leq M^{-1}(\ln3M)^{\frac{1}{16}c_0}.  \label{4.17}
\ee
In the following, we bound $b=v_re_r+v_ze_z$.\\
Denote $B_\rho(x_0)=\{x:|x-x_0|<\rho\}$, where $\rho>0$ to be determined later. By Biot-Savart law, $b$ satisfies
\be
-\Delta b={\rm curl}(w_\th e_\th).   \nn
\ee
From the estimates  of elliptic equation\cite{Gilbarg01}, for $q>1$,
\be
\sup\limits_{B_\rho(x_0)}|b|\leq C\left(\rho^{-\frac{3}{q}}\|b\|_{L^q(B_{2\rho}(x_0))}+\rho\sup\limits_{B_{2\rho}(x_0)}|w_\th|\right). \label{e4.12}
\ee
Take
\be
\rho=M^{-1}(\ln3M)^{\frac{1}{32}c_0},\  x_0\in\{(r,\th,z):r<\rho\}\  {\rm and}\ 1<q<2. \nn
\ee
By the assumption \eqref{e1.2} on $b$,
\bea
\rho^{-\frac{3}{q}}\|b\|_{L^q(B_{2\rho}(x_0))}&\leq&\rho^{-\frac{3}{q}}\|\frac{(\ln|\ln\frac{r}{3}|)^\al}{r}\|_{L^q(B_{2\rho}(x_0))} \nn  \\
&\leq& C \rho^{-\frac{3}{q}}\left[\int^{z_0+2\rho}_{z_0-2\rho}dz\int^{3\rho}_0\frac{(\ln|\ln\frac{r}{3}|)^{\al q}}{r^q}rdr\right]^{\frac{1}{q}} \nn  \\
&\leq& C \rho^{-\frac{2}{q}}\left[\int^{3\rho}_0\frac{(\ln|\ln\frac{r}{3}|)^{\al q}}{r^{q-1}}dr\right]^{\frac{1}{q}}.  \label{4.20}
\eea
We compute $\int^{3\rho}_0\frac{(\ln|\ln\frac{r}{3}|)^{\al q}}{r^{q-1}}dr$ as follow,
\bea
\int^{3\rho}_0\frac{(\ln|\ln\frac{r}{3}|)^{\al q}}{r^{q-1}}dr&=&\int^{+\infty}_{\frac{1}{3\rho}}(\ln\ln3r)^{\al q}r^{q-3}dr  \q r\ \ {\rm replaced}\ \ {\rm by}\ \ \f{1}{r}\nn \\
&=&\left(\int^{\frac{1}{3}e^{e^{\left(\ln\ln\frac{3}{\rho}\right)^2}}}_{\frac{1}{3\rho}}+\int^{+\infty}_{\frac{1}{3}e^{e^{\left(\ln\ln\frac{3}{\rho}\right)^2}}}\right)(\ln\ln3r)^{\al q}r^{q-3}dr       \nn  \\
&\leq&\left(\ln\ln\frac{3}{\rho}\right)^{2\al q}\int^{+\infty}_{\frac{1}{3\rho}}r^{q-3}dr+\int^{+\infty}_{\frac{1}{3}e^{e^{\left(\ln\ln\frac{3}{\rho}\right)^2}}}
\left(\frac{\ln\ln 3r}{r}\right)^{\al q}r^{q-3+\al q}dr.  \nn
\eea
Here $\frac{\ln\ln 3r}{r}$ is a decreasing function in the integral domain. Also we can pick a $q\in(1,2)$ such that $q-3+\al q<-1$. So
\bea
\int^{3\rho}_0\frac{\left(\ln|\ln\frac{r}{3}|\right)^{\al q}}{r^{q-1}}dr&\leq&C\left(\ln\ln\frac{3}{\rho}\right)^{2\al q}\left(\frac{1}{\rho}\right)^{q-2}+
\frac{(\ln\ln\frac{3}{\rho})^{2\al q}}{\left(\frac{1}{3}e^{e^{\left(\ln\ln\frac{3}{\rho}\right)^2}}\right)^{\al q}}\left(\frac{1}{3}e^{e^{\left(\ln\ln\frac{3}{\rho}\right)^2}}\right)^{q-2+\al q}
\nn  \\
&\leq&C\left(\ln\ln\frac{3}{\rho}\right)^{2\al q}\left(\frac{1}{\rho}\right)^{q-2}.  \label{4.21}
\eea
The last inequality holds because  $e^{e^{\left(\ln\ln\frac{3}{\rho}\right)^2}}>\frac{3}{\rho}$.\\
Since $\rho=M^{-1}(\ln3M)^{\frac{1}{32}c_0}$ and $\rho^{-1}\leq M$, so, from \eqref{4.20} and \eqref{4.21},
\bea
\rho^{-\frac{3}{q}}\|b\|_{L^q(B_{2\rho}(x_0))}&\leq& C \rho^{-\frac{2}{q}}\left(\ln\ln\frac{3}{\rho}\right)^{2\al}\left(\frac{1}{\rho}\right)^{1-\frac{2}{q}} \nn  \\
                                              &\leq&C\left(\ln\ln\frac{3}{\rho}\right)^{2\al}\left(\frac{1}{\rho}\right)  \nn \\
                                              &\leq&CM(\ln3M)^{-\frac{c_0}{32}}\cdot(\ln\ln3M)^{2\al}   \nn  \\
                                              &\leq&CM(\ln3M)^{-\frac{c_0}{64}}.  \label{e4.13}
\eea
While, due to \eqref{4.17}, when $x_0\in B_\rho(x_0)$,
\bea
\rho\sup\limits_{B_\rho(x_0)}|w_\th|&\leq& M^{-1}(\ln3M)^{\frac{c_0}{32}}M^2(\ln3M)^{-\frac{c_0}{16}} \nn  \\
       &\leq&CM(\ln3M)^{-\frac{c_0}{32}}. \label{e4.14}
\eea
Combining \eqref{e4.12},\eqref{e4.13} and \eqref{e4.14}, we have
\be
|b|\leq CM(\ln3M)^{-\frac{c_0}{64}} \q {\rm for}\ r<M^{-1}(\ln3M)^{\frac{1}{32}c_0};\label{e4.15}
\ee
next, when $M^{-1}(\ln3M)^{\frac{c_0}{32}}\leq r<1$,
\bea
&|b|\leq C\frac{(\ln\ln\frac{3}{r})^\al}{r}\leq CM(\ln3M)^{-\frac{c_0}{32}}(\ln\ln3M)^\al\leq CM(\ln3M)^{-\frac{c_0}{64}};  \label{4.22}
\eea
thirdly, when $1\leq r$,
\be
|b|\leq \f{C}{r}\leq C. \label{4.23}
\ee
Combining \eqref{e4.15},\eqref{4.22} and \eqref{4.23}, we get, for any $r>0$,
\be
|b|\leq CM(\ln3M)^{-\frac{c_0}{64}} . \label{4.24}
\ee
In the following,we bound $v_\th$.\\
Recall that $v_\th$ satisfies \eqref{e4.3}, then
\be
v_\th(r,z)=M|v^M_\th(rM,zM)|\leq M\left\{
\begin{aligned}
\min\left\{\frac{(\ln\frac{3}{r})^{-c_0}}{rM},rM\right\} \q&\q{\rm for}\ r<\frac{1}{M};  \\
\min\left\{\frac{(\ln\frac{3}{r})^{-c_0}}{rM},1\right\} \ \q&\q{\rm for}\ \frac{1}{M}<r<1.
\end{aligned}
\right.\nn
\ee
Firstly, when $r<\frac{1}{M}$, $\frac{(\ln\frac{3}{r})^{-c_0}}{rM}$ is a decreasing function while $rM$ is an increasing function with respect to $r$. Let $r_0$ be such that
\be
\frac{(\ln\frac{3}{r_0})^{-c_0}}{r_0M}=r_0M.  \nn
\ee
This gives
\be
(\frac{1}{r_0M})^{\frac{2}{c_0}}=\ln\frac{3}{r_0}>\ln 3M.    \nn
\ee
So
\be
r_0<(\ln3M)^{-\frac{c_0}{2}}M^{-1}.   \nn
\ee
Then
\be
|v_\th|\lesssim r_0M^2\lesssim(\ln3M)^{-\frac{c_0}{2}}M. \label{e4.16}
\ee
Next, when $\frac{1}{M}\leq r<1$,
\be
|v_\th|\lesssim\frac{(\ln\frac{3}{r})^{-c_0}}{r}\lesssim(\ln3M)^{-c_0}M. \label{e4.17}
\ee
Thirdly, When $1\leq r$,
\be
|v_\th|\leq \f{C}{r}\leq C. \label{e4.23}
\ee
Combining \eqref{e4.16},\eqref{e4.17} and \eqref{e4.23}, we have, for any $r>0$,
\be
|v_\th|\leq CM(\ln3M)^{-\frac{c_0}{64}}. \label{4.28}
\ee
Since $M$ is the maximum of $|v|$, $M=\max\{\sup |b|, \sup|v_\th|\}$. Due to the estimates \eqref{4.24} and \eqref{4.28}, we get
\be
M\leq CM(\ln3M)^{-\frac{c_0}{64}}.  \nn
\ee
This gives an upper bound for $M$ which completes the proof of Theorem 1.1.  \ef

\section{Proof of Theorem 1.3}

In this section ,we  prove Theorem 1.3 and give the estimate \eqref{e4.6} of the fundamental solution.

Following Davies \cite{Davies01} and Carlen-Loss \cite{Carlen01}, for a fixed constant vector $\al\in\bR^3$, let $\psi(x)=\al\cdot x$. For any $f\in C^\i_0(\bR^3;(0,+\i))$, define
\bea
P^\psi_{t,s}f(x)&=&e^{-\psi(x)}\int f(y)p(x,t;y,s)e^{\psi(y)}dy \nn \\
            &=&e^{-\al x}\int f(y)p(x,t;y,s)e^{\al y}dy \nn \\
            &\triangleq&f_t(x).  \nn
\eea

In fact, let $Q_R=B_R(0)\times(s,+\i)$ and define the Dirichlet fundamental solution $p^R(x,t;y,s)$ in $Q_R$ the same as {\bf Definition1.1} which satisfies the boundary condition
\be
p^R(x,t;y,s)|_{(x,t)\in \p B_R\times(s,+\i)}=0.  \nn
\ee
Due to the maximum principle, we have
\be
p^{R_1}(x,t;y,s)\leq p^{R_2}(x,t;y,s)\leq p(x,t;y,s),\q {\rm when}\ R_1<R_2. \nn
\ee
Also
\be
\lim\limits_{R\rightarrow +\i} p^R(x,t;y,s)=p(x,t;y,s).   \q a.e.  \nn
\ee

In a rigorous computation, all the integrals in the following should be done  in $B_R(0)$ with the function $f_t(x)$ replaced by $
f^R_t(x)\triangleq e^{-\al x}\int_{\bR^3} f(y)p^R(x,t;y,s)e^{\al y}dy $
which satisfies $f^R_t|_{x\in \p B_R}=0$ for all $t\geq s$. Then let $R\rightarrow +\i$ to reach the estimate of $p(x,t;y,s)$.  But for simplicity, we just carry out this process on $f_t(x)$ and assume that $f_t(x)$ vanishes on the boundary which means, $f_t(x)|_{|x|=+\i}=0$.\\

\noindent We divide the proof into 3 parts.\\
\noindent{\bf Part one: $L^2\rightarrow L^\i$ estimate of $P^\psi_{t,s}$}.

Let $k(t):[s,T]\rightarrow [2,\i]$ be a continuously differentiable increasing function to be determined later. By direct computation,we have
\bea
&&k(t)^2\pa f_t\pa^{k(t)-1}_{k(t)}\f{d}{dt}\pa f_t\pa_{k(t)}  \nn  \\
&=&k'(t)\int f^{k(t)}_t\ln\left(\f{f^{k(t)}_t}{\pa f_t\pa^{k(t)}_{k(t)}}\right)dx+k(t)^2\int f^{k(t)-1}_t\f{d}{dt}f_tdx \nn  \\
&=&k'(t)\int f^{k(t)}_t\ln\left(\f{f^{k(t)}_t}{\pa f_t\pa^{k(t)}_{k(t)}}\right)dx+k(t)^2\int(f^{k(t)-1}_te^{-\al x})(\Dl+\f{2}{r}\p_r-b\cdot\na)(e^{\al x}f_t)dx \nn
\\
&=&k'(t)\int f^{k(t)}_t\ln\left(\f{f^{k(t)}_t}{\pa f_t\pa^{k(t)}_{k(t)}}\right)dx+k(t)^2\left\{-\int\na(f^{k(t)-1}_te^{-\al x})\cdot\na(e^{\al x}f_t)\right. \nn  \\
&&\left.+\int\f{2}{r}f^{k(t)-1}_t\left[\p_r f_t+\p_r(\al x)f_t\right]-\int f^{k(t)-1}_t(b\cdot\al f_t+b\cdot \na f_t)\right\}  \nn  \\
&\triangleq&k'(t)\int f^{k(t)}_t\ln\left(\f{f^{k(t)}_t}{\pa f_t\pa^{k(t)}_{k(t)}}\right)dx+k(t)^2\{I+II+III\}.  \label{5.100}
\eea
Using Cauchy-Schwartz inequality, we have
\bea
I&=&-\int[(k(t)-1)f^{k(t)-2}_t\na f_t-f^{k(t)-1}_t\al]\cdot\left[\na f_t+\al f_t\right] \nn \\
&=&-\int\f{4(k(t)-1)}{k(t)^2}|\na f^{\f{k(t)}{2}}_t|^2-\int(k(t)-2)f^{k(t)-1}_t\al\cdot\na f_t+\int\al^2f^{k(t)}_t \nn  \\
&=&-\int\f{4(k(t)-1)}{k(t)^2}|\na f^{\f{k(t)}{2}}_t|^2-\int\f{2(k(t)-2)}{k(t)}f^{\f{k(t)}{2}}_t\al\cdot\na f^{\f{k(t)}{2}}_t+\int\al^2f^{k(t)}_t \nn \\
&\leq&-\int\f{4(k(t)-1)}{k(t)^2}|\na f^{\f{k(t)}{2}}_t|^2+\int\f{(k(t)-2)^2\ve}{k(t)^2}|\na f^{\f{k(t)}{2}}_t|^2+\int(1+\f{1}{\ve})\al^2f^{k(t)}_t, \label{ea.2}
\eea
where $\ve>0$ is to be determined later on.\\
Also,
\bea
II&=&\int\f{2}{r}f^{k(t)-1}_t[\p_r f_t+\p_r(\al x)f_t] \nn  \\
&\leq&\int\f{2}{k(t)}\p_rf^{k(t)}_tdrd\th dz+\int|\al|\f{2}{r}f^{k(t)}_t \nn \\
&=&-\f{2}{k(t)}\int dzd\th f^{k(t)}_t|_{r=0}+\int|\al|\f{2}{r}f^{k(t)}_t  \nn  \\
&\leq&\int|\al|\f{2}{r}f^{k(t)}_t. \label{ea.3}
\eea
The estimate  \eqref{ea.3} is due to our choice of $f\geq 0$, So $f_t|_{r=0}\geq 0$.\\
Moreover, due to the divergence-free property of $b$,
\bea
III&=&-\int f^{k(t)-1}_t(b\cdot\al f_t+b\cdot f_t)\} \nn  \\
&=&-\int f^{k(t)}_tb\cdot\al \nn \\
&\leq&\int (C_0+\f{1}{r})|\al|f^{k(t)}_t, \label{ea.4}
\eea
while using integration by parts,
\bea
\int \f{|\al|}{r}f^{k(t)}_t&=&|\int^{+\i}_{-\i}\int^{2\pi}_{0}\int^{+\i}_0 |\al|f^{k(t)}_tdrd\th dz| \nn  \\
&=&|-\int^{+\i}_{-\i}\int^{2\pi}_{0}\int^{+\i}_0\int2|\al|rf^{\f{k(t)}{2}}_t\p_rf^{\f{k(t)}{2}}_tdrd\th dz| \nn  \\
&\leq&\int2|\al|f^{\f{k(t)}{2}}_t\p_rf^{\f{k(t)}{2}}_tdx  \nn  \\
&\leq&\ve\int|\na f^{\f{k(t)}{2}}_t|^2+\f{\al^2}{\ve}\int f^{k(t)}_t. \label{ea.5}
\eea
Thus combining \eqref{ea.3},\eqref{ea.4} and \eqref{ea.5}, we have
\be
II+III\leq 3\ve\int|\na f^{\f{k(t)}{2}}_t|^2+3\f{\al^2}{\ve}\int f^{k(t)}_t+C_0|\al|\int f^{k(t)}_t. \label{5.50}
\ee
Now we recall the 3-d $log-Sobolev$ inequality.\\
\textit{ For all functions $u$ on $\bR^3$, together with their distributional gradients $\na u$ are square integrable, then}
\be
\int u^2\ln\left(\f{u^2}{\pa u\pa^2_2}\right)+\left(3+\f{3}{2}\ln a\right)\int u^2 \leq \f{a}{\pi}\int|\na u|^2  \label{5.5}
\ee
\textit{for all $a>0$.}\\
Inserting \eqref{ea.2},\eqref{5.50} and \eqref{5.5} into \eqref{5.100}, one has
\bea
k(t)^2\pa f_t\pa^{k(t)-1}_{k(t)}\f{d}{dt}\pa f_t\pa_{k(t)}&\leq& k'(t)\left[\f{a}{\pi}\int|\na f^{\f{k(t)}{2}}_t|^2-\left(3+\f{3}{2}\ln a\right)\int f^{k(t)}_t\right] \nn \\
&&+k(t)^2\left\{\left(\f{-4(k(t)-1)+\ve(k(t)-2)^2}{k(t)^2}+3\ve\right)\int|\na f^{\f{k(t)}{2}}_t|^2  \right.   \nn\\
&&\left.+\left((1+\f{4}{\ve})\al^2+C_0|\al|\right)\int f^{k(t)}_t\right\},  \nn
\eea
then
\bea
\pa f_t\pa^{k(t)-1}_{k(t)}\f{d}{dt}\pa f_t\pa_{k(t)}&\leq&\left[\f{k'(t)}{k(t)^2}\f{a}{\pi}-\f{4(k(t)-1)}{k(t)^2}+\f{(k(t)-2)^2}{k(t)^2}\ve+3\ve\right]\int|\na f^{\f{k(t)}{2}}_t|^2 \nn  \\
&&+\left[(1+\f{4}{\ve})\al^2+C_0|\al|-\f{k'(t)}{k(t)^2}\left(3+\f{3}{2}\ln a\right)\right]\int f^{k(t)}_t.  \nn
\eea

Here we can not choose $k(t):[s,T]\rightarrow[1,+\i]$ as Carlen-Loss in \cite{Carlen01} do. Because when $k(s)=1$, the coefficient of $\int|\na f^{\f{k(t)}{2}}_t|^2$ is $k'(s)\f{a}{\pi}+4\ve$ which is obviously positive when $k(t)$ is a continuously differentiable increasing function and $a>0$. It can not reach zero as Carlen-Loss in \cite{Carlen01} do. So we choose $k(t):[s,T]\rightarrow[2,+\i]$ to ensure the coefficient of $\int|\na f^{\f{k(t)}{2}}_t|^2$ is zero. \\
When $k(t)\in[2,\i)$,
\be
\f{k(t)-1}{k(t)}\geq \f{1}{2},\q \f{(k(t)-2)^2}{k(t)^2}<1.  \nn
\ee
So
\bea
\pa f_t\pa^{k(t)-1}_{k(t)}\f{d}{dt}\pa f_t\pa_{k(t)}&\leq&\left[\f{k'(t)}{k(t)^2}\f{a}{\pi}-\f{2}{k(t)}+4\ve\right]\int|\na f^{\f{k(t)}{2}}_t|^2  \nn  \\
&&+\left[(1+\f{4}{\ve})\al^2+C_0|\al|-\f{k'(t)}{k(t)^2}\left(3+\f{3}{2}\ln a\right)\right]\int f^{k(t)}_t.  \nn
\eea
Let $k(t)=2\s{\f{T}{T+s-t}}$ , $4\ve=\f{k'(t)}{k(t)^2}\f{a}{\pi}=\f{1}{k(t)}$, then
\be
a=\f{\pi k(t)}{k'(t)}=2\pi(T+s-t), \nn
\ee
\be
-\f{k'(t)}{k(t)^2}=-\f{1}{4\s{T(T+s-t)}}. \nn
\ee
Then
\bea
\pa f_t\pa^{-1}_{k(t)}\f{d}{dt}\pa f_t\pa_{k(t)}&\leq&\left[\left(1+32\s{\f{T}{T+s-t}}\right)\al^2+C_0|\al|\right.  \nn \\
                                                   &&\left.-\f{3}{4\s{T(T+s-t)}}-\f{3\ln(2\pi(T+s-t))}{8\s{T(T+s-t)}}\right]. \nn
\eea
Integrating the above inequality in $[s,T]$, we get
\bea
\ln\pa f_{T}\pa_{\i}-\ln\pa f_s\pa_{2}&\leq&\int^{T}_s\left[\left(1+32\s{\f{T}{T+s-t}}\right)\al^2+C_0|\al|\right.  \nn\\
                                    \hspace{50mm}&&-\f{3}{4\s{T(T+s-t)}}\left.-\f{3\ln(2\pi(T+s-t))}{8\s{T(T+s-t)}}\right]dt \nn \\
                                    &\leq&(\al^2+C_0|\al|+64\al^2)(T-s)-\f{3\ln 2\pi}{4}-\f{3\ln (T-s)}{4} \nn \\
                                    &\leq&(C_0|\al|+65\al^2)(T-s)+\ln(2\pi (T-s))^{-\f{3}{4}}. \nn
\eea
So
\be
\pa f_T\pa_{\i}\leq(2\pi (T-s))^{-\f{3}{4}}\exp\{(65\al^2+C_0|\al|)(T-s)\}\pa f\pa_{L^2}. \nn
\ee
That is
\be
\pa P^\psi_{t,s}f \pa_{\i}\leq(2\pi (t-s))^{-\f{3}{4}}\exp\{(65\al^2+C_0|\al|)(t-s)\}\pa f\pa_{L^2}.  \label{e5.50}
\ee\\
{ \bf Part two: $L^2\rightarrow L^\i$ estimate of the adjoint $( P^\psi_{t,s})^*$ of $P^\psi_{t,s}$}.

Now we come to investigate the adjoint $( P^\psi_{t,s})^*$ of $P^\psi_{t,s}$, for any $f,g\in\ C^\i_0(\bR^3)$,
\bea
(P^\psi_{t,s} f(x),g(x))&=&\int g(x)e^{-\psi(x)}\int f(y)e^{\psi(y)}p(x,t;y,s)dydx  \nn \\
&=&\int f(y)e^{\psi(y)}dy\int g(x)e^{-\psi(x)}p(x,t;y,s)dx \nn  \\
&\triangleq&((P^\psi_{t,s})^* g(y),f(y))  \nn
\eea
So
\be
(P^\psi_{t,s})^*g(y)=e^{\psi(y)}\int g(x)e^{-\psi(x)}p(x,t;y,s)dx. \nn
\ee
Here, note that we do not require $t\geq s$.
We denote $y=(y_1,y_2,y_3)$ and $y'=(y_1,y_2,0)$.

Let $p(x,t;y,s)$ be the fundamental solution of \eqref{1.9}, that is
\be
\p_t p(x,t;y,s)=\Delta_x p(x,t;y,s)-b\cdot\nabla_x p(x,t;y,s)+\frac{2}{r_x}\p_{r_x} p(x,t;y,s), \nn
\ee
when $t>s$. Here $r_x=\s{x^2_1+x^2_2}$. Then $p(x,t;y,s)$ ,with respect to $(y,s)$, satisfies
\be
-\p_s p(x,t;y,s)=\Delta_y p(x,t;y,s)+b\cdot\nabla_y p(x,t;y,s)-\frac{2}{r_y}\p_{r_y} p(x,t;y,s).  \nn
\ee
Let $\rho=-t,\tau=-s$. $p(x,\rho;y,\tau)$, with respect to $(y,\tau)$, satisfies
\be
\p_\tau p(x,\rho;y,\tau)=\Delta_y p(x,\rho;y,\tau)+b\cdot\nabla_y p(x,\rho;y,\tau)-\frac{2}{r_y}\p_{r_y} p(x,\rho;y,\tau).  \nn
\ee
 Let $p^*(y,\tau;x,\rho)=p(x,\rho;y,\tau)$, then
 \be
(P^\psi_{\rho,\tau})^*g(y)=e^{\psi(y)}\int g(x)e^{-\psi(x)}p^*(y,\tau;x,\rho)dx. \nn
 \ee
When $\tau>\rho$, $p^*(y,\tau;x,\rho)$ satisfies
\be
\p_\tau p^*(y,\tau;x,\rho)=\Delta_y p^*(y,\tau;x,\rho)+b\cdot\nabla_y p^*(y,\tau;x,\rho)-\frac{2}{r_y}\p_{r_y} p^*(y,\tau;x,\rho).  \nn
\ee
Then $p^*(y,\tau;x,\rho)$ is a fundamental solution of
\be
\p_\tau v=\Delta v+b\cdot\nabla v-\frac{2}{r}\p_{r} v,  \label{e5.70}
\ee
with respect to variables $(y,\tau)$ and $e^{-\psi(y)}(P^\psi_{\rho,\tau})^*g(y)$ is a solution of \eqref{e5.70}.

We now restrict the solution $v$ of \eqref{e5.70} such that $v(y,\tau)|_{|y'|=0}=0$. The reason is the following: let $v=rh$, then by direct computation, $h$ satisfies
\be
\p_s h=\Dl h-\f{1}{r^2}h+b\cdot\na h+\f{b_r}{r}h, \label{ea.9}
\ee
where $b=b_re_r+b_\th e_\th+b_ze_z$.

 If $|b|\leq C_0+\f{1}{r}$ , using Nash-Moser iteration argument as in the section 2 and noting that the integral of $\f{b_r}{r}h$ can be absorbed by that of $-\f{1}{r^2}h$ which is a good term in the energy estimate due to its minus sign. We can derive that the weak solution of \eqref{ea.9} is bounded.
So we can assume
\be
v|_{|y'|=0}=rh|_{|y'|=0}=0.   \nn
\ee
Then we have $(P^\psi_{t,s})^*g(y)|_{|y'|=0}=0$ when $s\geq t$.

Now we can follow the proof of $L^2\rightarrow L^\i$ estimate for $P^\psi_{t,s}$  to derive the $L^2\rightarrow L^\i$ estimate for $(P^\psi_{t,s})^*$. $(P^\psi_{t,s})^*$ has nearly the same form as $P^\psi_{t,s}$, but the signs on the terms $\f{2}{r}\p_r$ and $b\cdot\na$ are reversed. This makes the estimate a little different when we deal with the term II. If we denote
 \be
(P^\psi_{t,s})^*g(y)=g_s(y) \nn
 \ee
then
\bea
II&=-&\int\f{2}{r}g^{k(t)-1}_s[\p_r g_s+\p_r(\al x)] \nn  \\
&\leq&-\int\f{2}{k(t)}\p_rg^{k(t)}_tdrd\th dz+\int|\al|\f{2}{r}g^{k(t)}_s \nn \\
&=&\f{2}{k(t)}\int dzd\th g^{k(t)}_s|_{r=0}+\int|\al|\f{2}{r}g^{k(t)}_s  \nn  \\
&=&\int|\al|\f{2}{r}g^{k(t)}_s. \nn
\eea

Due to the vanishing property of $(P^\psi_{t,s})^*g(y)$  at $|y'|=0$, we can also get the estimate \eqref{ea.3} for $(P^\psi_{t,s})^*g(y)$, so the $L^2\rightarrow L^\i$ estimate \eqref{e5.50} is also right to $(P^\psi_{t,s})^*g(y)$.
\be
\pa(P^\psi_{t,s})^*g\pa_{\i}\leq(2\pi |t-s|)^{-\f{3}{4}}\exp\{(65\al^2+C_0|\al|)|t-s|\}\pa g\pa_{L^2}. \nn
\ee\\
\noindent{\bf Part three: $L^1\rightarrow L^\i$ estimate of $P^\psi_{t,s}$}.

Using the duality, we have the $L^1\rightarrow L^2$ estimate of $P^\psi_{t,s}$.
\be
\pa(P^\psi_{t,s})f\pa_{L^2}\leq(2\pi (t-s))^{-\f{3}{4}}\exp\{(65\al^2+C_0|\al|)(t-s)\}\pa f\pa_{L^1}. \nn
\ee
So
\bea
\pa(P^\psi_{t,s})f\pa_{L^\i}=\pa( P^\psi_{t,\f{t+s}{2}} P^\psi_{\f{t+s}{2},s})f)\pa&\leq&(2\pi (t-\f{t+s}{2}))^{-\f{3}{4}}\exp\{(65\al^2+C_0|\al|)(t-\f{t+s}{2})\}\pa P^\psi_{\f{t+s}{2},s}f\pa_{L^2}   \nn \\
&\leq&(\pi (t-s))^{-\f{3}{4}}\exp\{(65\al^2+C_0|\al|)\f{t-s}{2}\}\pa P^\psi_{\f{t+s}{2},s}f\pa_{L^2}.  \nn \\
&\leq&(\pi (t-s))^{-\f{3}{2}}\exp\{(65\al^2+C_0|\al|)(t-s)\}\pa f\pa_{L^1}.  \nn
\eea
This is equivalent to
\bea
p(x,t;y,s)&\leq&(\pi (t-s))^{-\f{3}{2}}\exp\{(65\al^2+C_0|\al|)(t-s)\}\exp\{\al(x-y)\}   \nn \\
          &\leq&(\pi (t-s))^{-\f{3}{2}}\exp\{65\al^2(t-s)+C_0|\al|(t-s)+\al(x-y)\}.     \nn
\eea
Let $\al=-\f{1}{65(t-s)}\f{x-y}{|x-y|}\left[|x-y|-C_0(t-s)\right]_+$.  With this choice of $\al$, we have
\be
\al\cdot(x-y)+C_0|\al|(t-s)+65\al^2(t-s)= -\f{1}{65t}\left[|x-y|-C_0(t-s)\right]^2_+,   \nn
\ee
then
\be
p(x,t;y,s)\leq (\pi (t-s))^{-\f{3}{2}}\exp\left\{-\f{1}{65(t-s)}\left[|x-y|-C_0(t-s)\right]^2_+\right\}. \nn
\ee
This gives the estimate \eqref{1.10} of $p(x,t;y,s)$. \\
Moreover, when $t\geq s$,
\bea
\p_t\int_{\bR^3}p(x,t;y,s)dx&=&\int_{\bR^3}\p_t p(x,t;y,s)dx \nn  \\
                          &=&\int_{\bR^3}(\Dl_x+\f{2}{r_x}\p_{r_x}-b\cdot\na_x)p(x,t;y,s)dx  \nn  \\
                          &=&\int^{+\i}_{-\i}\int^{2\pi}_0\int^{+\i}_0\f{2}{r_x}\p_{r_x} p(x,t;y,s)r_xdr_xd\th dz  \nn  \\
                          &=&\int^{+\i}_{-\i}\int^{2\pi}_0-2p(x,t;y,s)|_{r_x=0}d\th dz   \nn \\
                          &\leq&0,  \nn
\eea
so
\be
\int_{\bR^3}p(x,t;y,s)dx\leq\int_{\bR^3}p(x,s;y,s)dx=1. \nn
\ee
Also, when $t\leq s$,
\bea
\p_s\int_{\bR^3}p(x,t;y,s)dy&=&\int_{\bR^3}\p_s p(x,t;y,s)dy \nn  \\
                          &=&\int_{\bR^3}(\Dl_y-\f{2}{r_y}\p_{r_y}+b\cdot\na_y)p(x,t;y,s)dy  \nn  \\
                          &=&\int^{+\i}_{-\i}\int^{2\pi}_0\int^{+\i}_0-\f{2}{r_y}\p_{r_y} p(x,t;y,s)r_ydr_yd\th dz  \nn  \\
                          &=&\int^{+\i}_{-\i}\int^{2\pi}_02p(x,t;y,s)|_{r_y=0}d\th dz   \nn \\
                          &=&0,  \nn
\eea
so
\be
\int_{\bR^3}p(x,t;y,s)dy=\int_{\bR^3}p(x,t;y,t)dy=1. \nn
\ee
Thus we complete the proof of Theorem1.3. \ef

\indent

{\bf Acknowledgement.} I want to express my heartfelt gratitude to my advisor, Professor Qi S. Zhang, for his helpful discussion and guidance about this subject. At the same time I would like to thank Professor Zhen Lei in Fudan University for his inspiring idea in this work. I also want to express my warm thanks to my co-advisor, Professor Yin Huicheng, in Nanjing University for his constant encouragement.

\
\\

\end{document}